\documentclass[12pt]{amsart}
\usepackage{graphicx}
\usepackage{epsfig}
\usepackage{amsmath}
\usepackage{amssymb}
\usepackage{caption}

\usepackage[text={425pt,650pt},centering]{geometry}
\usepackage{color}
\usepackage{amssymb}
\usepackage{graphicx}
\usepackage{tikz}

\newtheorem{theorem}{Theorem}

\newtheorem{lemma}[theorem]{Lemma}

\theoremstyle{definition}

\newtheorem{definition}[theorem]{Definition}

\numberwithin{equation}{section}
\numberwithin{figure}{section}
\numberwithin{theorem}{section}

\newcommand{\nit}{\noindent}
\newcommand{\Z}{\mathbb{Z}}
\newcommand{\N}{\mathbb{N}}
\newcommand{\R}{\mathbb{R}}
\newcommand{\Q}{\mathbb{Q}}

\newcommand{\E}{\mathbb{E}}
\renewcommand{\P}{\mathbb{P}}
\newcommand{\F}{\mathcal{F}}
\newcommand{\A}{\mathcal{A}}

\newcommand{\Rd}{\mathbb{R}^d}

\newcommand{\ep}{\varepsilon}

\DeclareMathOperator*{\esssup}{ess\,sup}

\DeclareMathOperator{\BUC}{BUC}

\renewcommand{\bar}{\overline}
\renewcommand{\tilde}{\widetilde}

\begin{document}

\title[Homogenization of non-convex Hamilton-Jacobi equations]{Stochastic homogenization of nonconvex Hamilton-Jacobi equations in one space dimension}

\begin{abstract}
We  prove stochastic homogenization for a general class of coercive, nonconvex Hamilton-Jacobi equations in one space dimension. Some properties of the  effective Hamiltonian arising in the nonconvex case are also discussed.
\end{abstract}

\author[S. N. Armstrong]{Scott N. Armstrong}
\address{Ceremade (UMR CNRS 7534), Universit\'e Paris-Dauphine, Paris, France}
\email{armstrong@ceremade.dauphine.fr}

\author[H. V. Tran]{Hung V. Tran}
\address{Department of Mathematics\\
The University of Chicago\\ 5734 S. University Avenue Chicago, Illinois 60637, USA}
\email{hung@math.uchicago.edu}

\author[Y. Yu]{Yifeng Yu}
\address{Department of Mathematics\\
University of California at Irvine, California 92697, USA}
\email{yyu1@math.uci.edu}

\date{\today}

\maketitle

\section{Introduction}

\subsection{Motivation and overview}
We study the coercive Hamilton-Jacobi equation
\begin{equation} \label{e.pde}
u^\ep_t + H(Du^\ep)+ V\left( \frac x\ep  \right) = 0 \quad \mbox{in} \ \R \times (0,\infty).
\end{equation}
The Hamiltonian $H:\R \to \R$ is a deterministic function which satisfies $H(p) \to +\infty$ as $|p| \to +\infty$. In particular, we do not assume~$H$ is convex. The potential $V$ is a bounded, stationary random field sampled by an ergodic probability measure. We prove that, in the limit as the length scale $\ep > 0$ of the correlations tends to zero, the solution $u^\ep$ of~\eqref{e.pde}, subject to an appropriate initial condition, converges to the solution $u$ of the effective, deterministic equation
\begin{equation} \label{e.pdehom}
u_t + \overline H(Du) = 0 \quad \mbox{in} \ \R \times (0,\infty). 
\end{equation}

\smallskip

The random homogenization of Hamilton-Jacobi equations   has received much attention in the last fifteen years. The first  results were due to Rezakhanlou and Tarver~\cite{RT} and Souganidis~\cite{S}, who independently proved qualitative results for general convex, first-order Hamilton-Jacobi equations in stationary-ergodic setting. Later, these results were extended to the viscous case by Kosygina, Rezakhanlou and Varadhan~\cite{KRV} and, independently, by Lions and Souganidis~\cite{LS2} as well as to equations with time-dependent coefficients by Kosygina and Varadhan~\cite{KV} and Schwab~\cite{Sch}. New proofs of these results based on the notion of intrinsic distance functions appeared later in Armstrong and Souganidis~\cite{AS2} for the first-order case and in Armstrong and Tran~\cite{AT1} in the viscous case. The latter allowed for quantitative results, which appeared in Armstrong, Cardaliaguet and Souganidis~\cite{ACS} (see also Matic and Nolen~\cite{MN}) and Armstrong and Cardaliaguet~\cite{AC}.

\smallskip

There are however few such results for equations which are not convex (or, at least not quasi-convex, see~\cite{AS2,DS1}) in the gradient variable-- a fact which has been highlighted as one of the prominent  open problems in the field. Essentially the only previous result for a genuinely non-convex equation is due to the authors~\cite{ATY2013}. In that paper, we proved that the equation
\begin{equation*}
u_t^\ep + \left( |Du^\ep|^2-1 \right)^2 + V\left( \frac x\ep \right) = 0 \quad \mbox{in} \ \Rd \times (0,\infty)
\end{equation*}
homogenizes for stationary-ergodic potentials in all space dimensions $d\geq 1$. Using some ideas from~\cite{ATY2013}, in this paper we prove, in~$d=1$, that the special nonconvex gradient profile in the latter equation may be replaced by a general coercive function. Although our arguments are confined to one space dimension, this is the first stochastic homogenization result for a general class of nonconvex Hamilton-Jacobi equations.

\smallskip

As we will see, the main difficulty is to analyze the precise shock structure of solutions of~\eqref{e.pde}, in particular with the way the potential interacts ``non-locally" with the bumps in the graph of the Hamiltonian~$H$. We eventually argue by induction, removing some bumps at a time until we are left with a quasi-convex equation or the situation in Section 3 (the oscillation of $V$ is larger than the global oscillation of all such bumps) where homogenization result is obtained straightly. 

\subsection{Precise statement of the main result}
The random potential is modeled by a probability measure on the set of all potentials. More precisely, let\begin{equation*} \label{}
\Omega:=\BUC(\R)
\end{equation*}
be the space of real-valued, bounded and uniformly continuous functions on~$\R$. We define~$\F$ to be the~$\sigma$-algebra on~$\Omega$ generated by pointwise evaluations:
\begin{equation*} \label{}
\F := \, \mbox{$\sigma$--algebra generated by the family of maps} \quad \left\{ V \mapsto V(x) \,:\,  x\in \R \right\}. 
\end{equation*}
The translation group action of $\R$ on $\Omega$ is denoted by $\{ T_y\}_{y\in \R}$ where $T_y:\Omega \to \Omega$ is defined by
\begin{equation*} \label{}
\left( T_y V\right)(x) := V(x+y). 
\end{equation*}
We consider a probability measure $\P$ on $(\Omega,\F)$ satisfying the followings: there exists $\overline m>0$ such that 
\begin{equation} \label{e.pub}
\P \left[ \esssup_{x\in \R} V(x)=0 \right] =1  \quad \mbox{and} \quad \P \left[ \esssup_{x\in \R} (-V(x))=\overline m \right] =1 
\end{equation}
for every $E \in \F$ and $y\in \R$,
\begin{equation} \label{e.pstat}
\P \left[ E \right] = \P \left[ T_yE \right] \quad \mbox{(stationarity)}
\end{equation}
and
\begin{equation} \label{e.perg}
\P \big[ \cap_{z\in \R} T_zE \big]  \in \{ 0,1 \} \quad \mbox{(ergodicity).}
\end{equation}
Assume that $H\in C(\R)$ and
\begin{equation}\label{e.Hcoer}
\lim_{|p| \to \infty} H(p)=+\infty.
\end{equation}
Notice that, by ergodicity, there is no loss in generality in~\eqref{e.pub} compared to the assumption that $\P \left[ \esssup_{x\in\Rd} |V(x)| <\infty \right] =1$. 

\smallskip

We now present the main result. Throughout, all differential equations and inequalities in this paper are to be interpreted in the viscosity sense (see~\cite{EBook}). Recall that, for each~$\ep > 0 $ and~$g\in \BUC(\R)$, there exists a unique solution~$u^\ep(\cdot,g)\in C(\R \times[0,\infty))$ of~\eqref{e.pde} in~$\R \times (0,\infty)$, subject to the initial condition $u^\ep(x,0,g) = g(x)$. 

\begin{theorem}
\label{t.main}
Assume \eqref{e.pub}--\eqref{e.Hcoer} hold. Then there exists $\overline H \in C(\R)$ satisfying
\begin{equation} \label{e.Hbarcoer}
\overline H(p) \rightarrow +\infty \quad \mbox{as} \ |p| \to \infty
\end{equation}
such that, if we denote, for each $g\in \BUC(\R)$, the unique solution of~\eqref{e.pdehom} subject to the initial condition $u(x,0) = g(x)$ by $u(x,t,g)$, then
\begin{equation*} \label{}
\P \left[ \forall g\in \BUC(\R), \ \forall k>0, \ \limsup_{\ep \to 0} \sup_{(x,t) \in B_{k} \times [0,k]} \left| u^\ep(x,t,g) - u(x,t,g) \right| = 0 \right] = 1. 
\end{equation*}
\end{theorem}

\smallskip

We highlight two key properties of $\overline H$ which play significant roles in our proof.  

\smallskip

\nit   {\it $\bullet$  Quasi-convexification of the effective Hamiltonian}. As we will show, the effective Hamiltonian $\overline H(p)$  becomes quasi-convex  when $\bar m$ is large (see Theorem  \ref{thm:oscV-large}). Thus when the oscillation of the potential is large, we may expect the effective Hamiltonian to be ``less non-convex" than~$H$.  Similar facts have also been noticed in \cite{Q} and  \cite{ATY2013}. 

\smallskip

\nit {\it $\bullet$  Existence and nonexistence of sublinear correctors.}  In the random setting,  a simple example due to  Lions and Souganidis \cite{LS1} shows that the cell problem might not have sublinear solutions. As we will see, our proof actually demonstrates that, in $d=1$, sublinear correctors  exist away from those flat pieces of $\overline H$ where $\overline H$ attains local extreme values.

\section{Preliminaries}

\begin{definition}\label{def:reg-homo}   
We say that the pair $(H, V)$  is \emph{regularly homogenizable} (with respect to~$\P$) at $p\in  \R$  if $(H,V)$ satisfies \eqref{e.pub}--\eqref{e.Hcoer}, and there exists  $\overline H(p)\in   \R$ such that for  any $R>0$,
\begin{equation}\label{p-conv}
\P \left[  \limsup_{\lambda\to 0}\max_{|y|\leq {R/ \lambda}}\left|\lambda v_{\lambda}(y, p)+\overline H(p)\right|=0\right]=1,
\end{equation}
where  $v_{\lambda}(\cdot, p)$ is the unique continuous bounded viscosity solution to
\begin{equation}\label{Cell-p}
\lambda v_{\lambda}+H(p+v_{\lambda}')+V(y)=0   \quad \text{in $\R$}.
\end{equation}
 $(H, V)$  is called \emph{regularly homogenizable} if $(H, V)$ is regularly homogenizable at $p$ for every $p\in  \R$. 
\end{definition}

The merit of this definition is that the conclusion of Theorem~\ref{t.main} holds if $(H,V)$ is regularly homogenizable (see for example~\cite[Lemma 7.1]{ACS}). Moreover, in view of~\cite[Lemma 5.1]{AS2}, the condition~\eqref{p-conv}  is equivalent to the following seemingly weaker convergence assertion:
\begin{equation}\label{p-reduce}
\P\left[\lim_{\lambda\to 0}\left|\lambda v_{\lambda}(0, p)+\overline H(p)\right|=0\right]=1.
\end{equation}
In order to obtain Theorem \ref{t.main}, it is enough to prove the following statement. 

\begin{theorem}\label{thm:HV}
Assume \eqref{e.pub}--\eqref{e.Hcoer} hold. Then $(H,V)$ is regularly homogenizable.
\end{theorem}

We next notice that, by comparison principle, the property of being regularly homogenizable  is stable under the supremum norm. The proof is easy and thus we omit it. 

\begin{lemma}\label{lem:stability}  
Assume that $(H_n,V_n)$ is regularly homogenizable at $p\in  \R$ for each $n\in \N$, and there exists $(H,V)$ such that
\[
\lim_{n\to \infty}\left( \|H_n-H\|_{L^\infty(\R)}+\|V_n-V\|_{L^\infty(\R \times \Omega)}\right)=0.
\]
Then $(H, V)$ is also regularly homogenizable at $p\in \R$ and 
\[
\overline H(p)=\lim_{n\to \infty}\overline {H}_n(p).
\]
\end{lemma}

\noindent By Lemma \ref{lem:stability} and Lemma \ref{A5} in the appendix, we may assume in addition to \eqref{e.pub}--\eqref{e.Hcoer} the following assumptions throughout the paper
\smallskip

\noindent (H1) $H:\R \to [0,\infty)$ is Lipschitz continuous, $\min_{\R}H=H(0)=0$ and
\[
\lim_{|p| \to \infty} H(p)=+\infty.
\]
\noindent (H2) There exist $L \geq 0$ and $p_1>p_2>\ldots>p_{2L}>p_{2L+1}=0$ such that
\begin{itemize}
\item[(i)] $H$ is strictly increasing in $[p_1,\infty)$ and $[p_{2k+1},p_{2k}]$ for $1\leq k \leq L$,
\item[(ii)] $H$ is strictly decreasing in $[p_{2k},p_{2k-1}]$ for $1\leq k \leq L$.
\item[(iii)] $H(p_1), H(p_2),\ldots, H(p_{2L+1})$ are distinct positive numbers.
\end{itemize}
\noindent (H3) There exist $\tilde L \geq 0$ and $\tilde p_1<\tilde p_2 <\ldots <\tilde p_{2\tilde L+1}=0$ such that
\begin{itemize}
\item[(i)] $H$ is strictly decreasing in $(-\infty,\tilde p_1]$ and $[\tilde p_{2k},\tilde p_{2k+1}]$ for $1\leq k \leq \tilde L$,
\item[(ii)] $H$ is strictly increasing in $[\tilde p_{2k-1},\tilde p_{2k}]$ for $1\leq k \leq \tilde L$.
\item[(iii)] $H(\tilde p_1), H(\tilde p_2),\ldots, H(\tilde p_{2\tilde L+1})$ are distinct positive numbers. 
\end{itemize}
\noindent (H4) $V\in C^\infty(\R)$ and each level set of $V$ has no cluster points, that is, there does not exist any $y\in \R$ such that $V^{(k)}(y)=0$ for all $k\in \N$.
\smallskip

Set
\[
m_i:=H(p_{2i-1}) \quad \text{and} \quad M_i:=H(p_{2i})\quad \text{for} \ i=1,\ldots,L.
\]
We denote
\[
\phi_1:=H|_{[p_1,\infty)}, \quad \phi_i:=H|_{[p_i,p_{i-1}]}\quad \text{for}\ 2\leq i \leq 2L+1,
\]
and
\begin{itemize}
\item $\psi_1:[m_1,\infty) \to [p_1,\infty)$ as $\psi_1:=\phi_1^{-1}$,
\item $\psi_{2i}:[m_i,M_i] \to [p_{2i},p_{2i-1}]$ as $\psi_{2i}:=\phi_{2i}^{-1}$ for $1\leq i \leq L$,
\item $\psi_{2i-1}:[m_i,M_{i-1}] \to [p_{2i-1},p_{2i-2}]$ as $\psi_{2i-1}:=\phi_{2i-1}^{-1}$ for $2\leq i \leq L+1$.
\end{itemize}

\begin{center}
\includegraphics[scale=0.7]{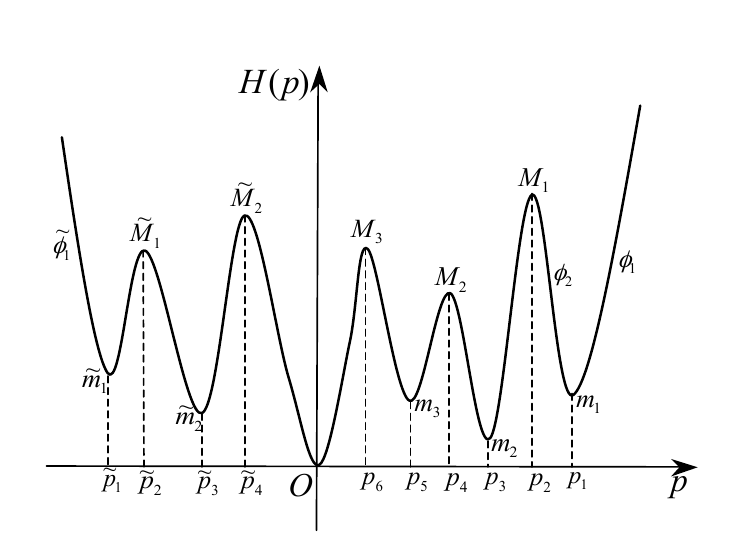}
\captionof{figure}{Graph of $H$ with $L=3$ and $\tilde L=2$}.
\end{center}

\begin{definition}\label{def:admissible}
We say that $f\in \mathcal A(H,V,\mu)$ for $\mu\in \R$ if $f\in L^\infty(\R)$ and any $u \in C^{0,1}(\R)$ solution to $u'=f$ is a solution to
\begin{equation}\label{metric-mu}
H(u')+V(y)=\mu \quad \text{in} \ \R.
\end{equation}
\end{definition}

\begin{lemma}\label{lem:comparison-R} 
Suppose that $u$ and $v$ are both viscosity solutions to
\[
\lambda w+H(p+w')+V(y)=0  \quad \text{in} \  B_{R/\lambda},
\]
for some $R>0$ and $\lambda \in (0,1)$.
Assume that $\lambda(|u|+|v|) \leq C$ in $B_{R/\lambda}$  and $\|H'\|_{L^\infty(\R)} \leq C$ for some $C>0$. Then
\begin{equation}\label{bound-uv}
\lambda |u(y)- v(y)|\leq   \frac{C}{R} \left(|y|^2+1\right)^{1/2}+\frac{C^2}{R} \quad \text{for} \ y\in B_{R/\lambda}.
\end{equation}
\end{lemma}

\begin{proof}
Let $w(y):=v(y)+\frac{C}{R} \left(|y|^2+1\right)^{1/2}+\frac{C^2}{R\lambda}$ for $y\in B_{R/\lambda}$.  
It is straightforward that $w$ is a viscosity supersolution to
\[
\lambda w+H(p+w')+V(y)=0   \quad \text{in} \  B_{R/\lambda},
\]
and $u\leq w$ on $\partial B_{R/ \lambda}$.  Thus $u \leq w$ in $B_{R/\lambda}$.
\end{proof}

\begin{lemma}[Generalized mean value theorem]\label{lem:meanvalue}
 Suppose that $u\in C([0,1],\R)$ and, for some $a,b \in \R$,
 \[
 u'(0^+)=\lim_{x \to 0^+} \frac{u(x)-u(0)}{x}=a \quad \text{and} \quad u'(1^-)=\lim_{x\to 1^-}\frac{u(1)-u(x)}{1-x}=b.
 \]
 Then:
 \smallskip
{\rm(i)} If $a<b$, then for any $c\in (a,b)$, there exists $x_c\in (0,1)$ such that $c \in D^- u(x_c)$, i.e., 
$$
u(x)\geq u(x_c)+c(x-x_c)-o(|x-x_c|) \quad \text{for} \ x\in (0,1).
$$
{\rm(ii)} If $a>b$, then for any $c\in (b,a)$, there exists $x_c\in (0,1)$ such that $c \in D^+ u(x_c)$, i.e., 
$$
u(x)\leq u(x_c)+c(x-x_c)+o(|x-x_c|)\quad \text{for} \ x\in (0,1).
$$
\end{lemma}

\begin{proof}
 It is enough to prove (i). For $c\in (a,b)$, set $w(x):=u(x)-cx$ for $x\in [0,1]$. There exists $x_c\in [0,1]$ such that
\[
 w(x_c)=\min_{x\in [0,1]} w(x).
\]
Note that $x_c \neq 0$ and $x_c\neq 1$ as $c\in (a,b)$. Thus $x_c \in (0,1)$, which of course yields that $c \in D^- u(x_c)$.
\end{proof}

\section{Homogenization in case the oscillation of $V$ is large}
In this section, we assume that $\tilde L=0$, and 
\begin{equation}\label{V-osc-large}
 \overline m > \max_{1\leq i,j \leq L} (M_i - m_j),
\end{equation}
and set $m_{\min}:=\min_{1\leq i \leq L} m_i>0$, $M_{\max}:=\max_{1\leq i \leq L} M_i>0$, and
\[
 \mathcal P:=\{\mu \geq 0\,:\,\mu \in (m_{\min}-\overline m, M_{\max})\}.
\]
\begin{definition}\label{def:decomposition}
 For $\mu \in \mathcal P$, a collection of finite intervals $\{I_i\}_{i\in \Z}$ is called a \emph{$(V,\mu)$-admissible decomposition} of $\R$ if
\begin{align*}
 &I_i=(a_i,a_{i+1}), \quad \lim_{i \to \pm \infty} a_i=\pm \infty,\quad \mu - V(a_i) \in \{m_j,M_j\,:\,1\leq j \leq L\},\\
 \quad &\text{and}\quad \{\mu-V(y)\,:\,y\in I_i\} \cap \{m_j,M_j\,:\,1\leq j \leq L\}=\emptyset.
\end{align*}
\end{definition}
\noindent Owing to \eqref{V-osc-large}, (H4), and Lemma \ref{A6}, $\{I_i\}$ exists and is unique up to a translation of indices in $\Z$. 
Furthermore, for any $i\in \Z$ and $y\in \R$,
\[
 T_y I_i= I_i+y.
\]

\begin{definition}\label{def:admissible-Vlarge}
 For $\mu \in \mathcal P$, we say $f\in \A(H,V,\mu)$ is furthermore \emph{ $(I_i,V,\mu)$-admissible} if 
\begin{align*}
& 0 \leq  f \leq \max\{p\geq 0\,:\, H(p) \leq \mu+\overline m\},\\
& f|_{I_i}=\psi_{j_i}(\mu-V) \quad \text{for some} \ j_i \in \{1,\ldots,2L+1\}, \ \text{for all} \ i \in \Z.
\end{align*}
\end{definition}

\begin{lemma}\label{lem:existence-P}
 For each $\mu \in \mathcal P$, there exists an  $(I_i,V,\mu)$-admissible function $f$.
\end{lemma}

\begin{proof}
 In view of Lemmas \ref{A1}, \ref{A2} in the appendix, there exist a strictly increasing sequence $\{b_i\}_{i\in \Z}$ and a Lipschitz continuous solution $u$ to
\eqref{metric-mu} such that $\lim_{i \to \pm \infty} b_i=\pm \infty$, $u\in C^1((b_i,b_{i+1}))$ for all $i\in \Z$ and
\[
 u'|_{(b_i,b_{i+1})}=\psi_{k_i}(\mu-V) \quad \text{for some} \ k_i \in \{1,\ldots,2L+1\}.
\]
By refinement, we may assume further that for each $i\in \Z$,
\[
 (b_i,b_{i+1}) \subseteq I_{l_i} \quad \text{for some} \ l_i \in \Z.
\]
For each $j\in \Z$, set
\[
 s_j:=\min\{k_i\,:\,(b_i,b_{i+1}) \subseteq I_j\}
\]
and
\[
 f=\psi_{s_j}(\mu-V) \quad \text{in} \ I_j.
\]
In light of one of the homotopy results, Lemma \ref{A3} in the appendix, we conclude that $f\in \A(H,V,\mu)$ and furthermore $(I_i,V,\mu)$-admissible.
\end{proof}

\noindent We now begin the identification of the effective Hamiltonian $\overline H$. For $\mu \in [0,\infty)\setminus \mathcal P$, we set
\begin{equation*}
 f_\mu:=\begin{cases}
         \psi_{2L+1}(\mu-V) \qquad &\text{if} \ \mu\leq  m_{\min}-\overline m ,\\
	 \psi_1(\mu-V) \qquad &\text{if} \ \mu \geq M_{\max}.
        \end{cases}
\end{equation*}
It is clear that $f_\mu \in \A(H,V,\mu)$ for $\mu \in [0,\infty)\setminus \mathcal P$.

\smallskip

For $\mu \in \mathcal P$ and $y\in \R$, define
\[
 \overline f_\mu(y):=\sup\{f(y)\,:\,f \ \text{is} \ (I_i,V,\mu)\text{-admissible}\}
\]
and
\[
 \underline f_\mu(y):=\inf\{f(y)\,:\,f \ \text{is} \ (I_i,V,\mu)\text{-admissible}\}.
\]

\begin{lemma}\label{lem:f-mu}
 Both $\overline f_\mu$ and $\underline f_\mu$ are stationary as well as $(I_i,V,\mu)$-admissible.
\end{lemma}

\begin{proof}
 Stationarity of $\overline f_\mu$ and $\underline f_\mu$ is straightforward. We now only check that $\overline f_\mu$ is $(I_i,V,\mu)$-admissible. We notice first that for all $i\in \Z$,
\[
 \overline f_\mu|_{I_i}=\psi_{j_i}(\mu-V) \quad \text{for some}\ j_i \in \{1,\ldots,2L+1\}.
\]
Thus, we only need to check that  for $u\in C^{0,1}(\R)$ such that $u'=\overline f_\mu$, $u$ is a solution of \eqref{metric-mu}  at $y=a_i$.

\smallskip

Pick $f_1, f_2$ which are $(I_i,V,\mu)$-admissible
such that 
\[
 \overline f_\mu=f_1 \quad \text{in} \ I_{i-1} \quad \text{and} \quad \overline f_\mu=f_2 \quad \text{in} \ I_i.
\]
{\it Case 1.} If
\[
 f_1(a_i^-):=\lim_{y\to a_i^-}f_1(y) \geq f_2(a_i^+):=\lim_{y\to a_i^+} f_2(y),
\]
then it is clear that
\[
 D^+u(a_i)=[\overline f_\mu(a_i^+),\overline f_\mu(a_i^-)] =[f_2(a_i^+),f_1(a_i^-)] \subseteq [f_1(a_i^+),f_1(a_i^-)].  
\]
{\it Case 2.} If
\[
 f_1(a_i^-):=\lim_{y\to a_i^-}f_1(y) \leq f_2(a_i^+):=\lim_{y\to a_i^+} f_2(y),
\]
then it is clear that
\[
 D^- u(a_i)=[\overline f_\mu(a_i^-),\overline f_\mu(a_i^+)] =[f_1(a_i^-),f_2(a_i^+)] \subseteq [f_2(a_i^-),f_2(a_i^+)].  
\]
The desired result follows.
\end{proof}

\begin{lemma}\label{lem:transition}
 For $\mu \in \mathcal P$ and $p\in \left[\E \left[  \underline f_\mu(0) \right],\E \left[  \overline f_\mu(0) \right]\right]$,
there exists a stationary function $f$ such that $f\in \A(H,V,\mu)$ and
\[
 p=\E \left[  f(0) \right] .
\]

\end{lemma}

\begin{proof}
 For $i\in \Z$, denote
\[
 \overline d_i:=\int_{a_i}^{a_{i+1}} \overline f_\mu(y)\,dy \quad \text{and} \quad
\underline d_i=\int_{a_i}^{a_{i+1}} \underline f_\mu(y)\,dy.
\]
According to \eqref{V-osc-large} and Lemma \ref{A6}, there exists a subsequence of intervals $\{I_{k_j}\}_{j\in \Z}$ 
such that $\lim_{j \to \pm \infty} k_j=\pm \infty$ and
\[
 \overline f_\mu=\underline f_\mu=\psi_1(\mu-V) \quad \text{in} \ I_{k_j}
\]
or
\[
 \overline f_\mu=\underline f_\mu=\psi_{2L+1}(\mu-V) \quad \text{in} \ I_{k_j}.
\]
By annexation if necessary, we may assume that for all $i \in \Z$
\[
 \begin{cases}
  \overline f_\mu=\underline f_\mu \qquad &\text{in} \ I_{2i}\\
  \overline f_\mu>\underline f_\mu \qquad &\text{in} \ I_{2i+1}.
 \end{cases}
\]
For $t\in [0,1]$ and $i\in \Z$, set $d_i(t):=t \overline d_i + (1-t) \underline d_i$ and
\[
 f_{\mu,t}:=\begin{cases}
       \overline f_\mu=\underline f_\mu \qquad &\text{in} \ I_{2i},\\
       f_{d_{2i+1}(t)}\left(\overline f_\mu,\underline f_\mu, I_{2i+1}\right)\qquad &\text{in} \ I_{2i+1}.
      \end{cases}
\]
By Lemma \ref{A4}, $f_{\mu,t} \in \A(H,V,\mu)$. The usual ergodic theorem gives
$$
\E\left[ \lim_{T\to \pm\infty}\frac1T\int_{0}^{T}f_{\mu,t}(y)\,dy=\mathbb {E}(f_{\mu,t}(0))\right]=1.
$$
So  it is  clear from the construction that the map $t \mapsto E(t):=\E \left[ f_{\mu,t}(0) \right]$
is Lipschitz continuous with 
\[
 E(0)=\E \left[  \underline f_\mu(0) \right] \quad \text{and} \quad E(1)=\E \left[  \overline f_\mu(0)\right],
\]
which gives us the desired result.
\end{proof}

The following lemma is similar to Lemma \ref{lem:f-mu}.

\begin{lemma}\label{lem:lim-stability}
 Assume that $\{\mu_m\}_{m\in \N}$ is a nonnegative sequence converging to $\mu$ and $f_m$ is $(I_i,V,\mu_m)$-admissible for each $m\in \N$.  Then we have that  

\smallskip

{\rm (1)} if $\mu\in  \mathcal {P}$,  then

\[
\limsup_{m\to \infty} f_m, \quad \liminf_{m\to\infty} f_m \quad \text{are} \quad (I_i,V,\mu) \text{-admissible};
\]

{\rm (2)}  if $m_{\min}\geq \bar m$ and $\mu\leq m_{\min}-\bar m$, then except on a countable set,
$$
\limsup_{m\to \infty}  f_m=\liminf_{m\to \infty}  f_m=\psi_{2L+1}(\mu-V);
$$

{\rm (3)}  if $\mu\geq M_{\max}$, then except on a countable set,
$$
\limsup_{m\to \infty}  f_m=\liminf_{m\to \infty}  f_m=\psi_{1}(\mu-V).
$$

\end{lemma}

\begin{proof}  

 Denote $f= \limsup_{m\to \infty}  f_m$.  The proof for $\liminf$ is similar.
 
 \smallskip

(1) Assume  $\mu\in  \mathcal {P}$.   Let $\{I_i\}_{-\infty<i<\infty}$  be the $(\mu,  V)$-admissible decomposition of $\R$.  For fixed $k\in \Z$ and $\ep>0$,  when $m$ is large enough, 
\[
\begin{cases}
\{\mu_m-V(y)\,:\, y\in (a_k+\ep, a_{k+1}-\ep)\} \cap \{M_i, m_i| \  1\leq i\leq L\}=\emptyset,\\
\{\mu_m-V(y)\,:\, y\in (a_{k+1}+\ep, a_{k+2}-\ep)\} \cap \{M_i, m_i| \  1\leq i\leq L\}=\emptyset.
\end{cases}
\]
Hence we can find four natural numbers $1\leq  l , \tilde l,  q , \tilde q\leq 2L+1$ and  two subsequences $\{f_{l_n}\}_{n\geq 1}$ and $\{f_{q_n}\}_{n\geq 1}$ such that
$$
f|_{I_{k}}=\psi_{l}(\mu-V)   \quad \mathrm{and}  \quad f|_{I_{k+1}}=\psi_{q}(\mu-V),
$$
$$
f_{l_n}=
\begin{cases}
\psi_{l}(\mu-V)   \quad \text{ in $(a_k+{1\over n},  a_{k+1}-{1\over n})$}\\
\psi_{\tilde l}(\mu-V)    \quad \text{ in $(a_{k+1}+{1\over n},  a_{k+2}-{1\over n})$}
\end{cases}
$$
and
$$
f_{q_n}=
\begin{cases}
\psi_{\tilde q}(\mu-V)   \quad \text{ in $(a_k+{1\over n},  a_{k+1}-{1\over n})$}\\
\psi_{q}(\mu-V)    \quad \text{ in $(a_{k+1}+{1\over n},  a_{k+2}-{1\over n})$}.
\end{cases}
$$
It suffices to show that for $u\in C^{0,1}(\R)$ such that $u'=f$, then $u$ is a solution of \eqref{metric-mu} at $a_{k+1}$.   Consider $u_{l}\in C^{0,1}(a_k, a_{k+2})$ with
$$
u_{l}':=
\begin{cases}
\psi_{l}(\mu-V)   \quad \text{in $I_k$}\\
\psi_{\tilde l}(\mu-V)  \quad \text{in $I_{k+1}$}
\end{cases}
$$
and
 $u_{q}\in C^{0,1}(a_k, a_{k+2})$ with
$$
u_{q}':=
\begin{cases}
\psi_{\tilde q}(\mu-V)   \quad \text{in $I_k$}\\
\psi_{q}(\mu-V)  \quad \text{in $I_{k+1}$}
\end{cases}
$$
Due to the stability of viscosity solutions,  $u_{l}$ and $u_{q}$ are both viscosity solutions to
\eqref{metric-mu} in $(a_k,a_{k+2})$.
By using the similar proof as the last part of that of Lemma \ref{lem:f-mu}, we are done.

\medskip

(2) Assume $m_{\min}\geq \bar m$ and $\mu\leq m_{min}-\bar m$.   Note if $u'\geq 0$ and $u$ is a solution of \eqref{metric-mu},  we must have that $u'=\psi_{2L+1}(\mu-V)$.  So it is clear that (2) holds for $x\in \R\backslash {A}$ where
$$
A:=\{y\in  \R\,:\,  V(y)=-\bar m\},
$$
which is either an empty set or a countable set due to (H4).

\medskip

(3) Assume  $\mu\geq M_{\max}$.   Note if $u'\geq 0$ and $u$ is a solution of \eqref{metric-mu},  we must have that $u'=\psi_{1}(\mu-V)$.  So it is clear that (3) holds for $y\in \R\backslash {B}$ where
$$
B:=\{y\in  \R\,:\,  V(y)=0\}
$$
which is either an empty set or a countable set due to (H4).
\end{proof}

For each $\mu \geq 0$, define
\[
 \mathcal I_\mu:=\begin{cases}
                  \left[\E \left[ \underline f_\mu(0) \right] , \E\left[ \overline f_\mu(0 \right] \right] \qquad &\text{for} \ \mu \in \mathcal P,\\ 
		  \vspace{1ex}
		  \left\{\E \left[   f_\mu(0)\right] \right\} \qquad &\text{for} \ \mu \in [0,\infty) \setminus \mathcal P.
                 \end{cases}
\]
For $\mu \in [0,\infty) \setminus \mathcal P$, we also write for consistency that
\[
 \E \left[  \underline f_\mu(0) \right]=\E \left[  \overline f_\mu(0) \right]=\E \left[   f_\mu(0) \right].
\]
Observe that  Lemma \ref{lem:transition} implies that
$$
\begin{array}{ll}
p\in \mathcal  I_{\mu}  &   \  \Rightarrow\ \text{existence of sublinear solutions to the cell problem}\\[3mm]
&\  \Rightarrow\  \text{$(H,V)$ is regularly homogenizable at $p$ and $\overline H(p)=\mu$}.
\end{array}
$$
In particular,  if $\mathcal  I_{\mu}$ is not a single point,  we obtain a flat piece. These intervals are mutually disjoint:

\begin{lemma}\label{lem:disjoint-I}
 If $\mu,\nu \in [0,\infty)$ with $\mu \neq \nu$, then $\mathcal I_\mu \cap \mathcal I_\nu=\emptyset$.
\end{lemma}

\begin{lemma}\label{lem:I-fill}
Set $q_0:=\E \left[  \underline f_0(0) \right]$. Then
\[
 \bigcup_{\mu \geq 0} \mathcal I_\mu=[q_0,\infty).
\]
\end{lemma}

\begin{proof} We divide the proof into two steps.

{\it Step 1.} We first show that  those intervals $\mathcal I_\mu$ are upper-semicontinuous with respect to $\mu$, i.e.,  for any nonnegative sequence $\{\mu_m\}$ converging to $\mu$
\begin{equation}\label{p1}
\begin{cases}
\E \left[\overline {f}_{\mu}(0) \right]\geq  \limsup_{m\to \infty}\E \left[\overline {f}_{\mu_m}(0) \right]\\
\E \left[\underline {f}_{\mu}(0) \right]\leq \liminf_{m\to \infty} \E \left[\underline {f}_{\mu_m}(0) \right].
\end{cases}
\end{equation}
In fact,  owing to Lemma \ref{lem:lim-stability},  it is obvious that 
$$
\underline {f}_{\mu} \leq \liminf_{m\to \infty} {\underline f}_{\mu_m} \leq \limsup_{m\to \infty} \overline f_{\mu_m}\leq \overline f_{\mu} \quad \text{a.e. in} \ \R.
$$
Hence using stationary ergodicity
\begin{align*}
\limsup_{m\to \infty}\E \left[\overline {f}_{\mu_m}(0) \right]&=\limsup_{m\to \infty}\int_{0}^{1}\E \left[\overline {f}_{\mu_m}(y) \right] dy\\
&\leq \int_{0}^{1}\E \left[\limsup_{m\to \infty}\overline {f}_{\mu_m}(y) \right] dy\\
&\leq \int_{0}^{1}\E \left[\overline {f}_{\mu}(y) \right] dy=\E \left[\overline {f}_{\mu}(0) \right].
\end{align*}
Similarly,  we can show that
$$
\liminf_{m\to \infty}\E \left[\underline {f}_{\mu_m}(0) \right]
\geq \E \left[\underline {f}_{\mu}(0) \right].
$$

{\it Step 2.}   This part  is similar to the proof of the intermediate value theorem for continuous functions.  We argue by contradiction.  If  the conclusion of this lemma is not true, then there exists $\overline p>q_0$ such that  $\overline p\notin \mathcal I_{\mu}$ for all $\mu\geq 0$.  

 For $\mu, \tilde \mu\geq 0$, if
$$
\max\{a\,:\, a\in \mathcal I_{\mu}\}< \overline p < \min \{a\,:\,a\in  \mathcal I_{{\tilde\mu}}\},
$$
then we compare $\overline p$ with the endpoints of $\mathcal I_{\mu+\tilde \mu\over 2}$.  By repeating this procedure,  we can find two sequences $\{\mu_n\}_{n\geq 1}$ and $\{\tilde {\mu}_{n}\}_{n\geq 1}$ such that
$$
\lim_{n\to \infty} \mu_n=\lim_{n\to \infty} {\tilde \mu_n}= \overline \mu\geq 0
$$
and for all $n\in \N$
$$
\max\{a\,:\, a\in \mathcal I_{\mu_n}\}< \overline p < \min \{a\,:\, a\in \mathcal  I_{{\tilde\mu}_n}\}.
$$
Then (\ref{p1})  implies that $\overline p\in  \mathcal I_{\overline \mu}$, which is a contradiction.
\end{proof}

We recall now that $\tilde L=0$ and thus $H$ is strictly decreasing on $(-\infty,0]$. Let 
\[
\Psi:=\left(H|_{(-\infty,0]}\right)^{-1} \quad \text{and} \quad q_{-1}:=\E \left[  \Psi(-V(0) \right]. 
\]
Sublinear correctors might not exist when $p\in  [q_{-1}, q_0]$.  Therefore,   we need to build a family of subsolutions which are sufficient to get the homogenization result at the minimum level $\{\overline H=0\}$.

\begin{lemma}\label{lem:subsln-0}
 For any $p\in [q_{-1},q_0]$ and $\delta>0$ sufficiently small, there exists a stationary function $f$ such that
\[
 p=\E \left[  f(0) \right]
\]
and for any $u\in C^{0,1}(\R)$ with $u'=f$, $u$ is a viscosity subsolution of
\begin{equation}\label{metric-delta}
 H(u')+V(y) =\delta \qquad \text{in} \ \R.
\end{equation}
\end{lemma}

\begin{proof}
 Choose $\delta$ such that
\[
 0<\delta <\frac{1}{2} \min\left\{\overline m, m_{\min}\right\},
\]
which implies that $\{p\,:\,H(p)<\delta\}$ is an interval containing $0$.

Take $\{b_i\}_{i\in \Z}$ to be a strictly increasing sequence satisfying $\lim_{i \to \pm \infty} b_i=\pm \infty$,
$V(b_i)=-\delta/4$, and
\[
 -\frac{\delta}{4} \notin \left\{V(y)\,:\,y\in (b_i,b_{i+1})\right\}.
\]
By (H4) and Lemma \ref{A6}, $\{b_i\}$ exists and is unique up to a translation of indices in $\Z$. 
For each $i\in \Z$, denote
\[
 \overline r_i:=\int_{b_i}^{b_{i+1}} \underline f_0(y)\,dy \quad \text{and} \quad
\underline r_i:=\int_{b_i}^{b_{i+1}} \Psi(-V(y))\,dy.
\]
For $t\in [0,1]$ and $i\in \Z$, set $r_i(t):=t \overline r_i+(1-t)\underline r_i$ and
\[
 f_{0,t}:=\begin{cases}
           t \underline f_0+(1-t)\Psi(-V) \quad &\text{in} \ (b_i,b_{i+1}) \ \text{if} \ V((b_i,b_{i+1}))\subset (-\delta/4,0],\\
           f_{r_i(t)}\left(\underline f_0,\Psi(-V),(b_i,b_{i+1})\right) \quad &\text{in} \ (b_i,b_{i+1}) \ \text{if} \ V((b_i,b_{i+1}))\subset (-\infty,-\delta/4).    
          \end{cases}
\]
Due to the choice of $\delta$,  we have that for any $u\in C^{0,1}(\R)$ such that $u'=f_{0,t}$, $u$ is a subsolution of \eqref{metric-delta}.
Repeating the last part of the proof of Lemma \ref{lem:transition}  yields the result.
\end{proof}

The following assertion holds in all dimensions $d\geq1$ provided~(H1) holds. 
\begin{lemma}\label{lem:lowerbound}
 Let $v_\lambda$ be the unique continuous bounded solution of \eqref{Cell-p} for some given $p\in \R$. Then
\[
 \P\left[\liminf_{\lambda \to 0} -\lambda v_\lambda(0,p) \geq 0\right]=1.
\]

\end{lemma}

\begin{theorem}\label{thm:oscV-large}
Assume  $\overline m \geq  \max_{1\leq i,j \leq L} (M_i - m_j)$. Then $(H,V)$ is regularly homogenizable and the effective Hamiltonian  $\overline H:\R\to [0,\infty)$  is quasi-convex.
\end{theorem}

\begin{proof}

Due to Lemma  \ref{lem:stability},  we may assume (\ref{V-osc-large}). 

\medskip

When  $p\geq \E \left[\underline {f}_{0}(0) \right]$, sublinear solutions  to the cell problem 
\begin{equation}\label{real-cell}
H(p+v')+V = \overline{H}(p) \quad \text{in} \ \R
\end{equation}
exist and are given by Lemma \ref{lem:transition}.

 When  $\E \left[\Psi(- V(0) \right]\leq p\leq \E \left[\underline {f}_{0}(0) \right]$,  sublinear solutions to cell problem \eqref{real-cell} might not exist.  However, combining  Lemma \ref{lem:subsln-0} and Lemma \ref{lem:lowerbound},   we have that $(H, V)$ is regularly homogenizable at $p$ and
$$
\overline H(p)=0
$$

  When  $p\leq \E \left[ \Psi(-V(0) \right]$,  $\overline H(p)\geq 0$ is the unique number given by 
$$
p=\E \left[  \Psi(\overline H(p)-V(0) )\right],
$$
and cell problem \eqref{real-cell}
has  a sublinear solution $v\in C^{0,1}(\R)$ with 
$$
v'=\Psi(\overline H(p)-V)-p \quad \text{in} \ \R.
$$

 It is clear  that such obtained $\overline H$ is quasiconvex:  $\overline H$ is increasing on  $[0,  \infty)$ and decreasing on $(-\infty, 0]$.
\end{proof}

\section{Homogenization by induction}
\subsection{Induction proof}
We present first the proof of Theorem \ref{thm:HV} by induction.
\begin{proof}[Proof of Theorem \ref{thm:HV}]
We prove by induction on $K:=\max\{L,\tilde L\}$.

\smallskip

{\it Base case.} If $K=0$,  the conclusion follows immediately  from Lemma \ref{lem:1d-AS2}, which is the one-dimensional case of~\cite{AS2}.

\smallskip

{\it Inductive hypothesis.} Assume that $(H,V)$ is regularly homogenizable for $K \leq k$ for some given $k \geq 0$. We now argue that $(H,V)$ is regularly homogenizable for $K=k+1$. Assume that $L \geq \tilde L$. In light of Lemma \ref{lem:gluing}, it suffices to show that $(H^+,V)$ is regularly homogenizable for
\[
H^+(p):=\begin{cases}
H(p) \qquad &\text{for} \ p \geq 0,\\
C|p| \qquad &\text{for} \ p \leq 0,
\end{cases}
\]
for some $C>\|H'\|_{L^\infty(\R)}$. There are two cases to be considered.

\smallskip

If $\overline m \geq \max_{1\leq i,j \leq L}(M_i-m_j)$, then the conclusion follows from Theorem \ref{thm:oscV-large}.
Otherwise, we use Lemmas \ref{lem:left} and \ref{lem:right} to reduce $H^+$ to simpler Hamiltonians and use the inductive hypothesis to achieve the result.
\end{proof}

\subsection{Gluing at the minimum point}
For some $C>\|H'\|_{L^\infty(\R)}$, define
\[
H^+(p):=\begin{cases}
H(p) \qquad &\text{for} \ p \geq 0,\\
C|p| \qquad &\text{for} \ p \leq 0,
\end{cases}
\]
and
\[
H^-(p):=\begin{cases}
H(p) \qquad &\text{for} \ p \leq 0,\\
C|p| \qquad &\text{for} \ p \geq 0.
\end{cases}
\]

\begin{lemma}\label{lem:gluing}
If both $(H^+,V)$ and $(H^-,V)$ are regularly homogenizable, then $(H,V)$ is also regularly homogenizable and moreover,
\begin{equation}\label{H-bar-glue}
\overline H(p)=\begin{cases}
\overline H^+(p) \qquad &\text{for}\ p \geq 0,\\
\overline H^-(p) \qquad &\text{for}\ p \leq 0.
\end{cases}
\end{equation}
Note that \eqref{H-bar-glue} is equivalent to the fact that $\overline H=\min\{\overline H^+,\overline H^-\}$.
\end{lemma}

\begin{proof}
It is enough to consider $p\geq 0$ and show that
\begin{equation}\label{conv-H-plus}
\P \left[ \lim_{\lambda \to 0}|\lambda v_\lambda(0,p)+\overline H^+(p)|=0\right]=1,
\end{equation}
where $v_\lambda$ is the solution of \eqref{Cell-p}. Let $v_\lambda^+$ be the solution of
\begin{equation}\label{Cell-plus}
\lambda v_\lambda^+ + H^+(p+(v_\lambda^+)')+V(y)=0 \quad \text{in} \ \R.
\end{equation}
By the usual comparison principle, 
\begin{equation}\label{gluing-0}
 \|\lambda v_\lambda(\cdot,p)\|_{L^\infty(\R)},  \|\lambda v_\lambda^+(\cdot,p)\|_{L^\infty(\R)} \leq H(p)+\overline m,
\end{equation}
and $v_\lambda \geq v_\lambda^+$ as $H^+ \geq H \geq 0$. Hence
\begin{equation}\label{gluing-2}
\P\left[\limsup_{\lambda \to 0} -\lambda v_{\lambda}(0,p) \leq \limsup_{\lambda \to 0} -\lambda v_\lambda^+ (0,p) = \overline H^+(p)\right]=1.
\end{equation}
When $\overline H^+(p)=0$,  \eqref{gluing-2} and Lemma \ref{lem:lowerbound} imply \eqref{conv-H-plus} immediately. Note also that
$\overline H^+(0)=0$. We thus only need to consider the case $p>0$ and $\overline H^+(p)>0$.

\smallskip

As $(H^+,V)$ is regularly homogenizable, 
\begin{equation}\label{gluing-3}
\P \left[ \forall R>0,\ \limsup_{\lambda \to 0} \max_{|y| \leq R/\lambda}\left|\lambda v_\lambda^+(y,p)+\overline H^+(p)\right|=0 \right] = 1.
\end{equation}
Fix $V$ belonging to this event. There exists $\lambda(R,V)>0$ such that for all $\lambda<\lambda(R,V)$ and $y\in B_{R/\lambda}$
\begin{equation}\label{gluing-4}
 -\lambda v_\lambda^+(y,p) \geq \frac{\overline H^+(p)}{2}.
\end{equation}
We claim that, for $R>8(H(p)+\overline m)/p$ and $\lambda <\lambda(R)$,
\begin{equation}\label{gluing-key}
 p+(v_\lambda^+(y,p))' \geq 0 \quad \text{for a.e.} \ y \in B_{R/(2\lambda)}.
\end{equation}
If \eqref{gluing-key} were false, there would exist $y_0 \in B_{R/(2\lambda)}$ such that 
$p+(v_\lambda^+(y,p))'(y_0) < 0$. On the other hand, by \eqref{gluing-0} and the choice of $R$,
\[
 \frac{2\lambda}{R} \int_{R/(2\lambda)}^{R/\lambda}\left(p+(v_\lambda^+)'(y)\right)\,dy \geq \frac{p}{2}>0.
\]
We use Lemma \ref{lem:meanvalue} to yield the existence of $y_1 \in (y_0,R/\lambda)$ such that $0 \in D^-u_{\lambda}^{+}(y_1)$ for $u_{\lambda}^{+}=py+v_\lambda^+(y,p)$,
and hence
\[
 \lambda v_\lambda^+(y_1,p)+H^+(0)+V(y_1) \geq 0,
\]
which contradicts \eqref{gluing-4}. Therefore, \eqref{gluing-key} is true and 
\[
 \lambda v_\lambda^+ + H(p+(v_\lambda^+)')+V(y)=0 \quad \text{in} \ B_{R/(2\lambda)}.
\]
The comparison result in Lemma \ref{lem:comparison-R} gives
\[
 \lambda\left|v_\lambda^+(0,p)-v_\lambda(0,p)\right| \leq \frac{C}{R}.
\]
Sending $R\to \infty$ to get the conclusion.
\end{proof}

\subsection{Gluing results in case the oscillation of $V$ is not large}
We assume that $\tilde L=0$ and
\begin{equation}\label{V-osc-small}
 \overline m < \max_{1\leq i,j \leq L} (M_i-m_j)
\end{equation}
in Lemma \ref{lem:left} and \ref{lem:right}. Due to (iii) in (H2),  there exist unique $k,L\in \{1,\ldots,L\}$ such that 
\[
 M_k=M_{\max} \quad \text{and} \quad m_l=m_{\min}.
\]
Then of course $\overline m < M_k-m_l=\max_{1\leq i,j \leq L} (M_i-m_j)$.
We need to consider two cases $l>k$ and $l\leq k$ as the nature of the difficulties is different.

\subsection{Left steep side}
We consider first the case that $l>k$.

\smallskip

Let $H_1:\R \to [m_l,\infty)$ be a coercive Lipschitz continuous function satisfying that $H_1\geq H$ and
\[
 \begin{cases}
  H_1=H \quad \text{on} \ (-\infty,p_{2k}],\\
  H_1 \quad \text{is strictly increasing in} \ (p_{2k},\infty),
 \end{cases}
\]
and $H_2:\R \to [m_l,\infty)$ be a coercive Lipschitz continuous function so that $H_2 \geq H$ and
\[
 \begin{cases}
  H_2=H \quad \text{on} \ [p_{2l},\infty),\\
  H_2 \quad \text{is strictly decreasing in} \ (-\infty,p_{2l}),
 \end{cases}
\]
and $H_3:=\max\{H_1,H_2\}$.

\begin{center}
\includegraphics[scale=0.7]{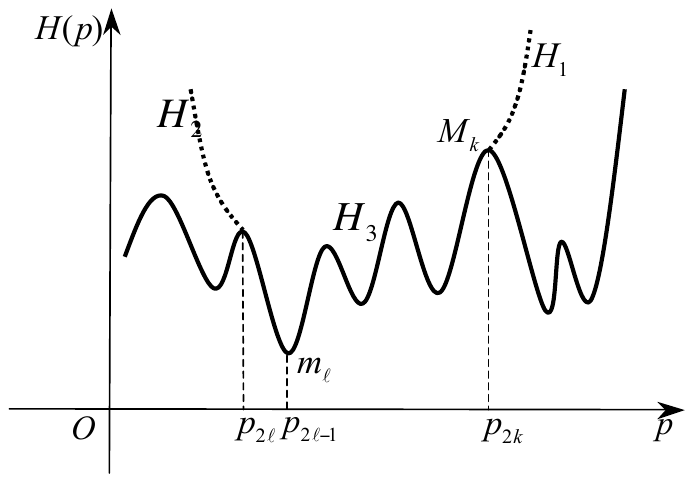}
\captionof{figure}{$H_1$, $H_2$ and $H_3$ in the gluing lemma \ref{lem:left}}.
\end{center}

\begin{lemma}\label{lem:left}
 Assume that $(H_i,V)$ are regularly homogenizable for all $i=1,2,3$. Then $(H,V)$ is also regularly homogenizable
and
\[
 \overline H=\min\{\overline H_1, \overline H_2\}.
\]
\end{lemma}

\begin{proof}     It is very easy to verify this lemma in  the periodic setting (i.e., $V$ is a periodic function).  To highlight the main ideas, we present it here.  

\medskip
\noindent {\bf Periodic case.} It is clear that
$$
\overline H\leq \min\{\overline H_1, \  \overline H_2\}.
$$
For fixed $p\in  \R$,  let $v(\cdot,  p)\in  C^{0,1}(\R)$ be a periodic viscosity solution to the cell problem
\eqref{real-cell}.
Then we must have either
$$
\{p+v'(y, p)\,:\,  y\in \R\} \subseteq (-\infty,  p_{2k}]
$$
or
$$
\{p+v'(y, p)\,:\,   y\in \R\} \subseteq [p_{2l},  \infty).
$$
Otherwise,  due to the periodicity and Lemma \ref{lem:meanvalue},  there exist $y_1, y_2\in  \R$ such that
$$
H(p_{2k})+V(y_1)\leq \overline H(p),  \quad  \mathrm{and}  \quad H(p_{2l-1})+V(y_2)\geq \overline H(p).
$$
So  $\overline m\geq M_{k}-m_l$, which is a contradiction.   Hence $v(\cdot,  p)$ is either a viscosity to
$$
H_1(p+v')+V(y)=\overline H(p)  \quad \text{in $\R$}
$$
or a viscosity solution to
$$
H_2(p+v')+V(y)=\overline H(p)   \quad \text{in $\R$}.
$$
Accordingly,  $\overline H(p)=\overline {H_1}(p)$ or $\overline H(p)=\overline {H_2}(p)$.

\medskip

\nit {\bf Random case.}    We note first that $\overline H_3(p_{2l-1})=\min \overline H_3=m_l$. Set
\[
 A:=\left\{p>p_{2l-1}\,:\,m_l<\overline H_3(p)<M_k-\overline m\right\}.
\]

{\it Step 1.} We first show that \eqref{p-reduce} holds for $p\in A$ and
\begin{equation}\label{left-1}
\overline H(p)= \overline H_i(p) \quad \text{for} \ i=1,2,3.
\end{equation}
The proof of this step is very similar to that of Lemma \ref{lem:gluing} hence is being sketched only.
As $(H_3,V)$ is regularly homogenizable, we have
\[
\P \left[ \forall R>0,\ \limsup_{\lambda \to 0} \max_{|y| \leq R/\lambda} |\lambda v_{3\lambda}(y,p)+\overline H_3(p)|=0 \right] = 1,
\]
where $v_{3\lambda}$ is the viscosity solution to
\begin{equation}\label{metric-H3}
\lambda v_{3\lambda}+H_3(p+v_{3\lambda}')+V(y)=0 \quad \text{in}\ \R.
\end{equation}
As usual, $\lambda \|v_{3\lambda}\|_{L^\infty(\R)} \leq H_3(p)+\overline m$. Set
\[
\delta:=\min\left\{\overline H_3(p)-m_l, M_k-\overline m - \overline H_3(p)\right\}.
\]
There exists $\lambda_3(R,V)>0$ such that when $\lambda<\lambda_3(R,V)$ 
\begin{equation}\label{left-2}
\max_{y\in B_{R/\lambda}} \left| \lambda v_{3\lambda}(y,p)+\overline H_3(p)\right| \leq \frac{\delta}{8}.
\end{equation}
Fix $R> 4(H_3(p)+\overline m)/(p-p_{2l-1})$ and $\lambda<\lambda_3(R,V)$. 
Then \eqref{left-2} yields, for $y\in B_{R/\lambda}$,
\[
H_3(p+v_{3\lambda}'(y,p)) \leq \overline H_3(p)+\frac{\delta}{8}+\overline m < M_k-\frac{\delta}{8}.
\]
Therefore, there exists $\tau>0$ depending only on $H,\delta$ such that
\begin{equation}\label{left-3}
p+v_{3\lambda}'(\cdot,p) <p_{2k}-\tau \quad \text{a.e. in} \ B_{R/\lambda}.
\end{equation}
On the other hand, the choice of $R$ allows us to get that
\[
\frac{2\lambda}{R}\int_{R/(2\lambda)}^{R/\lambda} \left(p+v_{3\lambda}'(y,p)\right)\,dy >p_{2l-1},
\]
which yields, by using the same proof as that of \eqref{gluing-key},
\begin{equation}\label{left-4}
p+v_{3\lambda}'(\cdot,p) >p_{2l-1} \quad \text{a.e. in} \ B_{R/(2\lambda)}.
\end{equation}
Combining \eqref{left-3} and \eqref{left-4} to achieve that
\begin{equation}\label{left-key-1}
p_{2l-1}<p+v_{3\lambda}'(\cdot,p) <p_{2k}-\tau \quad \text{a.e. in} \ B_{R/(2\lambda)},
\end{equation}
and thus
\[
\lambda v_{3\lambda}+H(p+v_{3\lambda}')+V(y)=0 \quad \text{in} \ B_{R/(2\lambda)}.
\]
So the comparison result in Lemma \ref{lem:comparison-R} yields
\[
\lambda \left|v_{3\lambda}(0,p)-v_\lambda(0,p)\right| \leq \frac{C}{R}.
\]
Letting $R\to \infty$ to conclude Step 1. Since $H_1=H_2=H_3$ in $[p_{2l},  p_{2k}]$,  (\ref{left-key-1}) immediately leads to $\overline {H_1}=\overline {H_2}=\overline {H_3}$ in $A$. 

\smallskip

{\it Step 2.} We claim that \eqref{p-reduce} holds for $p\leq p_{2l-1}$ and
\[
\overline H(p)=\overline H_1(p).
\]
This is due to $\overline m<M_k-m_l$.   Let $v_{1\lambda}$ be the unique viscosity solution to
\[
\lambda v_{1\lambda}+H_1(p+v_{1\lambda}')+V(y)=0 \quad \text{in} \ \R.
\]
Then 
\[
\P \left[ \forall R>0,\ \limsup_{\lambda \to 0} \max_{|y| \leq R/\lambda} |\lambda v_{1\lambda}(y,p)+\overline H_1(p)|=0 \right] = 1.
\]
Let $\delta_1:=M_k-m_l-\overline m>0$. For each $R>0$, there exists $\lambda_1(R,V)>0$ such that, for $\lambda<\lambda_1(R,V)$, 
\[
\max_{y\in B_{R/\delta}}\left|\lambda v_{1\lambda}(y,p)+\overline H_1(p)\right| \leq \frac{\delta_1}{8}.
\]
Choose $\tau_1>0$ sufficiently small so that
\[
H_1< m_l +\frac{\delta_1}{8} \qquad \text{in} \ (p_{2l-1}-\tau_1,p_{2l-1}+\tau_1).
\]
For $R>4(H_1(p)+\overline m)/\tau_1$ and $\lambda<\lambda_1(R,V)$, we also have the following key property
\begin{equation}\label{left-key-2}
p+v_{1\lambda}'(\cdot,p) \leq p_{2k} \quad \text{a.e. in} \ B_{R/(2\lambda)}.
\end{equation}
If not,  then there exists $y_0\in  B_{R/(2\lambda)}$ such that  $p+v_{1\lambda}'(y_0, p)>p_{2k}$.   Due to the choice of $R$,
$$
{2\lambda\over R}\int_{-{R/\lambda}}^{-{R/(2\lambda)}}v_{1\lambda}'(y, p)\,dy< \tau_1,  \quad {2\lambda\over R}\int_{R/(2\lambda)}^{R/\lambda}v_{1\lambda}'(y,p)\,dy<\tau_1.
$$
According to Lemma \ref{lem:meanvalue},  there must exists $y_{+}\in (y_0,  {R/\lambda})$ and $y_{-}\in (-{R/ \lambda}, y_0)$ such that
$$
H_1(p_{2l-1}+\tau_1)+V(y_{-})\geq -\lambda v_{1\lambda}(y_{-}, p)
$$
and
$$
H_1(p_{2k})+V(y_{+})\leq -\lambda v_{1\lambda}(y_{+}, p).
$$
Hence $\overline m\geq  M_k-m_l-{\delta_1/2}$, which contradicts the choice of $\delta_1$.  Thus, \eqref{left-key-2} holds and 
\[
\lambda v_{1\lambda}+H(p+v_{1\lambda}')+V(y)=0 \quad \text{in} \ B_{R/(2\lambda)},
\]
and, in light of Lemma \ref{lem:comparison-R},
\[
\lambda\left|v_{1\lambda}(0,p)-v_\lambda(0,p)\right| \leq \frac{C}{R}.
\]
Step 2 is complete.

\smallskip

{\it Step 3.} By similar arguments as in the above two steps, we can conclude that
\begin{itemize}
\item For $p \geq p_{2k}$ then \eqref{p-reduce} holds and $\overline H(p)=\overline H_2(p)$.
\item For $p\in \R$ such that $\overline H_1(p)<M_k-\overline m$, then \eqref{p-reduce} is true and $\overline H(p)=\overline H_1(p)$.
\item For $p \in [p_{2l-1},p_{2k}]$ with $\overline H_2(p)>m_l$, then \eqref{p-reduce} holds and $\overline H(p)=\overline H_2(p)$.
\item For $p\in [p_{2l-1},p_{2k}]$ and $\overline H_2(p)<M_k-\overline m$, then $\overline H_2(p)=\overline H_3(p)$.
\end{itemize}
Since $H_3=\max\{H_1,H_2\}$, one gets $\overline H_3 \geq \max\{\overline H_1,\overline H_2\}$.
In particular, if $p\in [p_{2l-1},p_{2k}]$ and $\overline H_3(p) \geq M_k -\overline m$, then by the last assertion above
\[
\overline H_3(p) \geq \overline H_2(p) \geq M_k - \overline m>m_l,
\]
and thus
\[
(p_{2l-1},p_{2k}) \subset A \cup \{p\,:\,\overline H_1(p)<M_k-\overline m\} \cup \{p\in [p_{2l-1},p_{2k}]\,:\,\overline H_2(p)>m_l\}.
\]
The proof is complete.
\end{proof}

\subsection{Right side is steeper} 
We consider now the case $l\leq k$.  We cannot simply copy the method when $l>k$.   The subtlety is that the decomposition in  the previous  case will not lead to a simpler Hamiltonian if $l=1$.  We need to  employ the following different approach.  

\smallskip

Let $H_1:\R \to [0,\infty)$ be a coercive Lipschitz continuous function satisfying that $H_1\geq H$ and
\[
 \begin{cases}
  H_1=H \quad \text{on} \ (-\infty,p_{2k}],\\
  H_1 \quad \text{is strictly increasing in} \ (p_{2k},\infty),
 \end{cases}
\]
and $H_2:\R \to [m_l,\infty)$ be a coercive Lipschitz continuous function so that $H_2 \geq H$ and
\[
 \begin{cases}
  H_2=H \quad \text{on} \ [p_{2k},\infty),\\
  H_2 \quad \text{is strictly decreasing in} \ (-\infty,p_{2k}).
 \end{cases}
\]

\begin{center}
\includegraphics[scale=0.7]{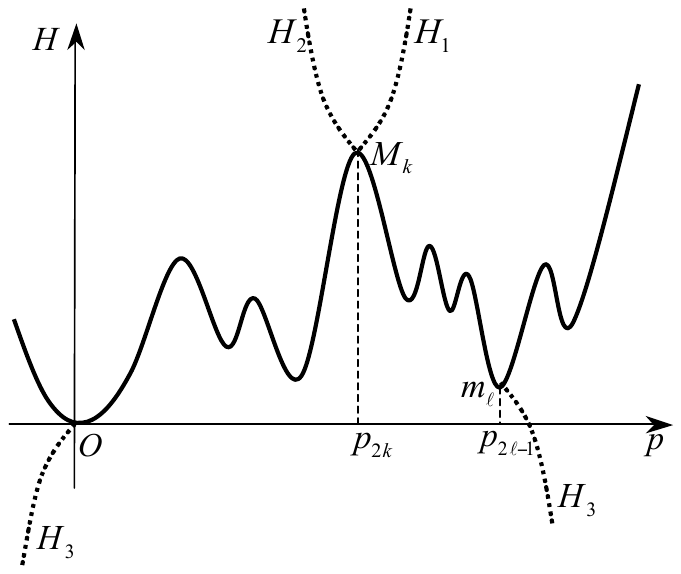}
\captionof{figure}{Graphs of $H_1$ and $H_2$ in Lemma \ref{lem:right}}.
\end{center}

\begin{lemma}\label{lem:right}
Assume both $(H_1,V)$ and $(H_2,V)$ are regularly homogenizable. Then $(H,V)$ is regularly homogenizable and
\[
\overline H(p)=\begin{cases}
\overline H_1(p) \qquad &\text{for} \ p\leq 0,\\
\min\{\overline H_1(p), \overline H_2(p),M_k-\overline m\} \qquad &\text{for} \ p \in [0,p_{2l-1}],\\
\overline H_2(p) \qquad &\text{for} \ p \geq p_{2l-1}.
\end{cases}
\]
\end{lemma} 

\begin{proof}  
As in the previous gluing lemma,   for readers' convenience,  we first prove the second equality in the statement when $V$ is periodic.
\medskip

\noindent{\bf  Periodic case.}   For $p\in  [0,  p_{2l-1}]$,   let $v(\cdot,  p)\in C^{0,1}(\R)$ be a periodic viscosity solution to \eqref{real-cell}.
Since 
\[
\int_{0}^{1}p+v'(y,p)\,dy=p\in  [0,  p_{2l-1}],
\]
  using Lemma \ref{lem:meanvalue}, it is easy to see that  $\overline H(p)\leq M_k-\overline m$ and therefore
$$
\overline H(p)\leq \min\{ \overline H_1(p), \overline H_2(p),M_k-\overline m\}.
$$
If $\overline H(p)<M_k-\overline m$,  then we have either
$$
\{p+v'(y, p)\,:\, y\in  \R\}\subseteq  (-\infty,  p_{2k})
$$
or
$$
\{p+v'(y, p)\,:\, y\in  \R\}\subseteq  (p_{2k},  \infty).
$$
Otherwise,  the periodicity and  Lemma \ref{lem:meanvalue}  imply the existence of $y_1\in  \R$ such that
$$
M_k-\overline m\leq H(p_{2k})+V(y_1)\leq \overline H(p),
$$
which contradicts  our  assumption.  Hence $v$ is either a viscosity solution to
$$
H_1(p+v')+V(y)=\overline H(p)  \quad \text{in  $\R$}
$$
or  $v$ is a viscosity solution to
$$
H_2(p+v')+V(y)=\overline H(p)  \quad \text{in  $\R$}.
$$
So  $\overline H(p)=\overline H_1(p)$ or $\overline H(p)=\overline H_2(p)$.

\medskip

\nit {\bf Random case}.  Proofs of the first and third  equalities  in the statement are similar to that of  Step 2  in  the proof of Lemma \ref{lem:left}.   We will prove the equality in the middle.    Using similar arguments to  Step 1  in the proof of Lemma \ref{lem:left},  we can deduce that

{\it Claim 1.} For $p\in  \R$,  if $\overline H_1(p)<M_k-\overline m$,  then $(H,V)$ is regularly homogenizable at $p$ and
$$
\overline H(p)=\overline H_1(p).
$$

{\it Claim 2.} For $p\in  \R$,  if $\overline H_2(p)<M_k-\overline m$,  then $(H,V)$ is regularly homogenizable at $p$ and
$$
\overline H(p)=\overline H_2(p).
$$
It is easy to see that $\overline H_1(0)=0$ and $\overline H_2(p_{2l-1})=m_l$.  Also, since
$$
\inf_{y\in  \R}\{H_1(p_{2k})+V(y)\}=\inf_{y\in  \R}\{H_2(p_{2k})+V(y)\}= M_k-\overline m,
$$
we have that  $\min\{\overline H_1(p_{2k}),\  \overline H_2(p_{2k})\}\geq M_k-\overline m$.  Now denote
$$
q_1:=\min \{p\in  [0,  p_{2k}]\,:\,  \overline H_1(p)=M_k-\overline m\}>0
$$
and
$$
q_2:=\max \{p\in  [p_{2k}, p_{2l-1}]\,:\,   \overline H_2(p)=M_k-\overline m\}<p_{2l-1}.
$$

Claim 1 and Claim 2 imply that  $(H, V)$ is regularly homogenizable for $p\in   [0,  q_1)\cup  (q_2,  p_{2l-1}]$ and
$$
\overline H(p)=
\begin{cases}
\overline H_1(p) \quad \text{when $p\in  [0, q_1)$}\\
\overline H_2(p)  \quad \text{when $p\in  (q_2,  p_{2l-1}]$}.
\end{cases}
$$
Our next goal is to show that $(H,V)$ is regularly homogenizable for $p\in  [q_1,  q_2]$ and
\begin{equation}\label{flat1}
\overline {H}|_{[q_1, q_2]}\equiv  M_k-\overline m.
\end{equation}

\medskip

Owing to Claims 1 and 2, and the stability Lemma \ref{lem:stability},  we have that $(H,V)$ are regularly homogenizable at $q_1$ and $q_2$ with
$$
\overline H(q_1)=\overline H(q_2)=M_k-\overline m.
$$
Now choose  $H_3:  \R\to  (-\infty,  M_k]$ to be Lipschitz continuous function such that $H\geq H_3$,  $\lim_{|p|\to +\infty} H_3(p)=-\infty$ and
$$
\begin{cases}
H_3=H  \quad \text{in $[0,  p_{2l-1}]$}\\
\text{ $H_3$ \quad is strictly increasing on $(-\infty, 0]$}\\
\text{ $H_3$ \quad is strictly decreasing on $[p_{2l-1},\infty)$}.
\end{cases}
$$
Using similar arguments as  Step 1 in the proof of Lemma \ref{lem:left}, we have that

\medskip

{\it Claim 3.} $(H_3, V)$ is regularly homogenizable at $q_1$ and $q_2$ and
$$
\overline H_3(q_1)=\overline H_3(q_2)=M_k-\overline m.
$$
Let $H_0(p):=-H_3(p_{2k}-p)+M_k$ for $p\in \R$.  It is easy to check that $w$ is the  viscosity solution to
$$
\lambda w+H_0(p+w')-V-\overline m=0  \quad \text{in $\R$}
$$
if and only if  $\tilde w=-w$ is a viscosity solution to
$$
\lambda \tilde w+H_3(p_{2k}-p+{\tilde w}')+V+\overline m-M_k=0  \quad \text{in $\R$}.
$$
By applying Lemma \ref{lem:r-flat}  to $(H_0, -V-\overline m)$,  we deduce that $(H_3, V)$ is regularly homogenizable at $p\in  [q_1, q_2]$ and
\begin{equation}\label{flatpart3}
\overline H_3|_{[q_1, q_2]}\equiv  M_k-\overline m.
\end{equation}
 Let $v_{\lambda}(\cdot, p)\in C^{0,1}(\R)$ be the unique bounded viscosity solution to \eqref{Cell-p}.
Since $H\geq H_3$,   (\ref{flatpart3}) says that  for $p\in  [q_1, q_2]$
\begin{equation}\label{1side}
\mathbb {P}\left[  \liminf_{\lambda\to 0}-\lambda v_{\lambda}(0, p)\geq M_k-\overline m\right]=1.
\end{equation}
Now  choose $\tilde H:\R\to \R$ to be a Lipschitz continuous function  satisfying that  $H\leq \tilde H$,
$$
\tilde H(p_{2k})=M_k,  \quad  \tilde H(0)=0,  \quad \tilde H(p_{2l-1})=m_l
$$
and $\tilde H|_{(-\infty,  0]}$ is strictly deceasing,  $\tilde H|_{[0, p_{2k}]}$ is strictly increasing,  $\tilde H|_{[p_{2k}, p_{2l-1}]}$ is strictly decreasing and $\tilde H|_{[p_{2l-1}, \infty)}$ is strictly increasing (see the figure  below.)
\begin{center}
\includegraphics[scale=0.7]{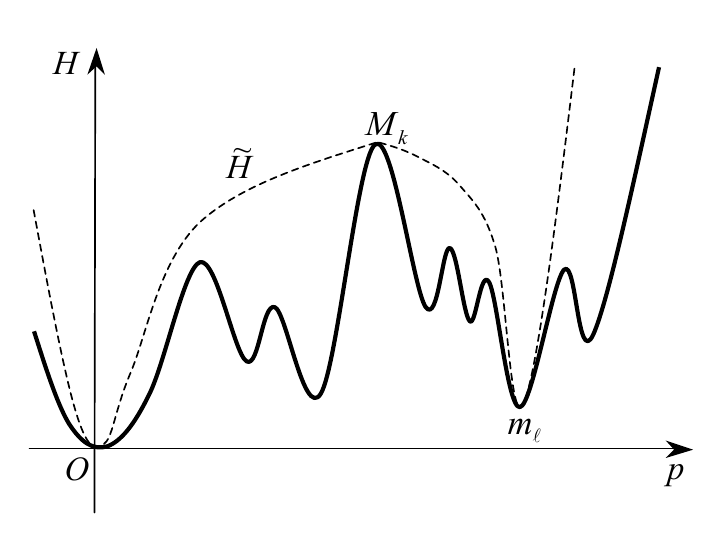}
\captionof{figure}{Graph of $\tilde H$}.
\end{center}
Since $\overline m<M_k-m_l<M_k$,  owing to Lemma \ref{localization},  $(\tilde H, V)$ is regularly homogenizable  and
$$
\overline {\tilde H}(p)\leq M_k-\overline m    \quad \text{for $p\in  [0,  p_{2l-1}]$}.
$$
Comparison principle implies that for $p\in  [0,  p_{2l-1}]$,
$$
\mathbb {P}\left[ \limsup_{\lambda\to 0}-\lambda v_{\lambda}(0, p)\leq M_k-\bar m\right]=1.
$$
Combining this with (\ref{1side}), we obtain~(\ref{flat1}). 
\end{proof}

\begin{lemma}\label{lem:r-flat}
Suppose that $(H,V)$ is regularly homogenizable at $q\in \R$ and $\overline H(q)=0$. Then $(H,V)$ is regularly homogenizable at all points $p\in I$, where $I$ is the line segment between $q$ and $0$, and
\[
\overline H(p)=0 \quad \text{for}\ p \in I.
\]
\end{lemma}

\begin{proof} As usual, we provide first the proof for the periodic case.
\medskip

\noindent {\bf Periodic Case.} When $V$ is periodic, the proof is quite simple.  Assume $q>0$.  It is obvious that $\overline H(p)\geq 0$.  So we only to verify that $\overline H(p)\leq 0$ for $p\in  [0, q]$.

Pick $y_0\in \R$ so that $V(y_0)=0=\min_{\R}V$.
Let $v(\cdot,  q)\in  C^{0,1}(\R)$ be a viscosity solution to cell problem
$$
H(q+v')+V(y)=0  \quad \text{in  $\R$}
$$
subject to the condition that $q y_0+v(y_0, q)=p-q$.  Then
\begin{equation}\label{uniformzero}
\lim_{r\to 0} \|q+v'(\cdot,q)\|_{L^{\infty}(B_r(y_0))}=0.
\end{equation}
For fixed $p\in  [0,  q]$,  set $w(y):=\max\{qy+v(y),  \  0\}$ in $[y_0,  y_0+1]$ and extend $w-py$ periodically to $\R$.  Note that \eqref{uniformzero} implies that $w$ is differentiable at $y_0$ and $w'(y_0)=0$.  Then $h=w-py$ is a  periodic viscosity subsolution to
$$
H(p+h')+V(y)=0  \quad \text{in  $\R$}.
$$
Thus $\overline H(p)\leq 0$.  

\medskip

\nit {\bf Random Case.} It is enough to consider the case where $q>0$.  Denote
$$
M^{+}=\max_{1\leq i\leq L}M_i  \quad \mathrm{and} \quad M^{-}=\max_{1\leq i\leq\tilde  L}\tilde M_i.
$$
If $\max\{M^{+}, \  M^{-}\} \leq \bar m$,  (1) follows immediately from  Theorem \ref{thm:oscV-large}.    Let us consider the case
$$
\min\{M^{+}, \  M^{-}\}> \bar m.
$$
The case that one of them is no larger than $\bar m$ is simpler.   Write
$$
k_+=\max\{1\leq i\leq L|\   M_i>\bar m\},   \quad  k_-=\max\{1\leq i\leq\tilde  L|\  \tilde  M_i>\bar m\}.
$$
Let $v_\lambda(\cdot,q)$ be the solution of \eqref{Cell-p} with $p=q$.
By the hypothesis,
\[
\P \left[ \forall R>0,\ \lim_{\lambda \to 0} \sup_{y\in B_{R/\lambda}}|\lambda v_\lambda(y,q)|=0 \right] = 1
\]
and by the ergodic theorem
\begin{equation}\label{r-flat-ergodic}
\lim_{s \to  \pm\infty} \frac{1}{s} \int_0^s \mathbf{1}_{\{y\,:\,-\delta<V(y)\leq 0\}}\,dy=\E \left[\left| \{y\,:\,-\delta<V(y)\leq 0\} \right| \right]>0,
\end{equation}
where
\[
\delta =\frac{1}{4}\min \left\{M_{k_{+}}-\overline m,\ \tilde  M_{k_{-}}-\overline m, \ \overline m,\  m_{\min}\right\}.
\]
There exists $\lambda(R,V)>0$ such that for $\lambda<\lambda(R, V)$
\[
 \sup_{y\in B_{R/\lambda}}|\lambda v_\lambda(y,q)| <\delta.
\]
In view of \eqref{r-flat-ergodic}, we can choose a sequence $\{\lambda_n\} \to 0$ such that for all $n\in \N$, $\lambda_n \in (0,\lambda(\delta,R))$ and
\[
\int_{R/(2\lambda_n)}^{R/\lambda_n} \mathbf{1}_{\{y\,:\,-\delta<V(y)\leq 0\}}\,dy,  \quad   \int_{-{R/\lambda_n}}^{-{R/(2\lambda_n)}} \mathbf{1}_{\{y\,:\,-\delta<V(y)\leq 0\}}\,dy>0.
\] 
Pick $y_{1n}^{+}\in (R/(2\lambda_n), R/\lambda_n)$ and $y_{1n}^{-}\in (-R/\lambda_n, -R/2(\lambda_n))$such that $v_\lambda(\cdot,q)$ is differentiable at $y_{1n}^{\pm}$ and $V(y_{1n}^{\pm}) \in (-\delta,0)$. Therefore, $H(q+v_{\lambda_n}^{'}(y_{1n}^{\pm})) \leq 2\delta$ and 
\begin{equation}\label{r-flat-1}
q+v_{\lambda_n}'(y_{1n}^{\pm})\in  (\tilde p_{2\tilde L}, p_{2L}).
\end{equation}
On the other hand, for all $y\in B_{R/\lambda_n}$, one has
\begin{equation}\label{r-flat-2}
H(q+v_{\lambda_n}'(y)) \leq \delta + \overline m \leq \min\{ M_{k^{+}}-3 \delta,  \tilde  M_{k_{-}}-3\delta\}.
\end{equation}
We combine \eqref{r-flat-1}, \eqref{r-flat-2}, and Lemma \ref{lem:meanvalue} to deduce that 
\begin{equation}\label{r-flat-3}
q+v_{\lambda_n}'(y) \in (\tilde p_{2k_{-}},\ p_{2k_{+}}) \quad \text{a.e. in} \ B_{R/(2\lambda_n)}.
\end{equation}

Let $\hat H:\R \to [0,\infty)$ be a Lipschitz continuous function satisfying that $\hat H \geq H$ and
$$
\hat H=H  \quad \text{on $[\tilde p_{2k_{-}},  p_{2k_{+}}]$},
$$
$\hat H|_{[p_{2k_{+}}, \infty)}$ is strictly increasing and $\tilde H|_{(-\infty,  \tilde p_{2k_{-}}]}$ is strictly decreasing.  See the following figure. 
\begin{center}
\includegraphics[scale=0.7]{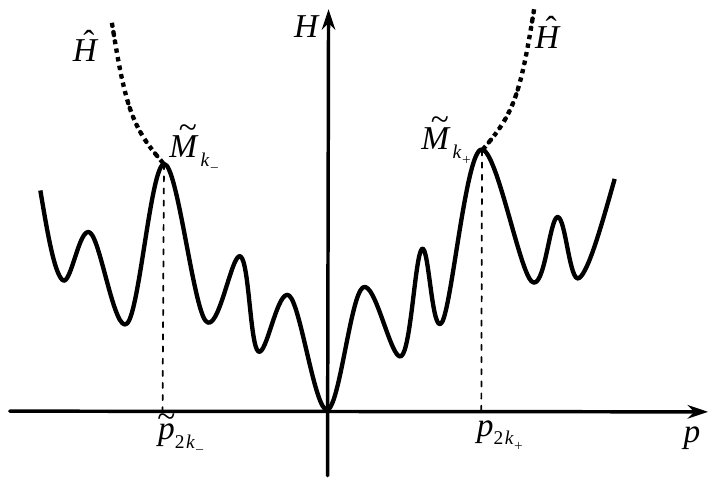}
\captionof{figure}{Graph of $\hat H$}.
\end{center}
By Theorem \ref{thm:oscV-large} and Lemma \ref{lem:gluing}, $(\hat H,V)$ is regularly homogenizable and $\overline {\hat H}:\R\to  [0,\infty)$ is quasi-convex. 
Let $\hat v_{\lambda_n}$ be the solution of
\begin{equation}\label{metric-hat}
\lambda_n \hat v_{\lambda_n}(y,q)+\hat H(q+\hat v_{\lambda_n}')+V(y)=0 \quad \text{in}\ \R.
\end{equation}
Note that both $\hat v_{\lambda_n} (\cdot,q)$ and $v_{\lambda_n}(\cdot,q)$ are solutions of \eqref{metric-hat} in $B_{R/(2\lambda_n)}$ by \eqref{r-flat-3}. We apply Lemma \ref{lem:comparison-R} to yield
\[
\lambda_n |v_{\lambda_n}(0,q)-\hat v_{\lambda_n}(0,q)| \leq \frac{C}{R},
\]
which of course gives us that $\overline{\hat H}(q)=0$. Thus, $\overline{\hat H}=0$ on $[0,q]$. Since  $\hat H \geq H$, we have that
$$
\mathbb{P}\left[\limsup_{\lambda\to 0}-\lambda v_{\lambda}(0,p)\leq 0\right]=1.
$$
Combining with Lemma \ref{lem:lowerbound}, the conclusion follows. 
\end{proof}

\section{Explicit formula of $\overline H$ in case the oscillation of $V$ is small}
The following lemma is the 1d case of   \cite{AS2}.  Since the proof is very easy,  we present it here for completeness. 

\begin{lemma}\label{lem:1d-AS2}
Suppose that $L=\tilde L=0$.  Then $(H,V)$ is regularly homogenizable and  the formula of $\overline H$ is given as follows
\begin{equation}
\begin{cases}
p=\E \left[\psi_{1}(\overline H(p)-V(0)) \right]  \quad \text{for $p\geq \E\left[\psi_1(-V(0))\right]$}\\
\overline H(p)\equiv 0    \quad \text{for $p\in [\E\left[\Psi(-V(0))\right], \E \left[\psi_1(-V(0))\right] ]$}\\
p=\E \left[\Psi(\overline H(p)-V(0)) \right]  \quad \text{for $p\leq \mathbb {E}\left[\Psi(-V(0))\right]$}.
\end{cases}
\end{equation}
Here $\Psi=H^{-1}:[0,\infty)\to  (-\infty,  0]$. 
\end{lemma}

\begin{proof}  We only need to prove the middle equality since the other two are obvious due to the existence of sublinear correctors.  For $t\in [0,1]$,  denote $u(t):=tu_{+}+(1-t)u_{-}$  where 
$$
u_{+}(y):=\int_{0}^{y}\psi_1(-V(z))\,dz  \quad \mathrm{and}  \quad u_{-}(y):=\int_{0}^{y}\Psi(-V(z))\,dz.
$$
Clearly,  $u(t)$ is a viscosity subsolution to
$$
H(u(t)')+V(y)=0  \quad \text{in $\R$}.
$$
Moreover,  $u(t)'=t\psi_1(-V)+(1-t)\Psi(-V)$ is stationary and $\mathbb {E}\left[u(t)'(0)\right]=p(t)$ where
\[
p(t)=t\E \left[\psi_1(-V(0))\right]+(1-t)E\left[\Psi(-V(0))\right].
\]
  So we have that
$$
\mathbb{P}\left[\limsup_{\lambda\to 0}-\lambda v_{\lambda}(0,p(t))\leq 0\right]=1.
$$
Combining this with  Lemma \ref{lem:lowerbound} yields the middle equality. 
\end{proof}

Now let us look at the case   $\tilde L=0$  and  $L\geq 1$.  For convenience,  set  $m_{L+1}=0$.  
We assume in this section that 
\begin{equation}\label{explicit-small}
\overline m < \min_{1\leq k \leq L} \min \{M_k-m_k, M_k-m_{k+1}\}.
\end{equation}
We denote, for $1\leq k \leq L$,
\begin{align*}
&p_{2k-1}^{+}:=\E \left[\psi_{2k-1}(m_k-V(0)) \right],\quad p_{2k-1}^{-}:=\E \left[\psi_{2k}(m_k-V(0)) \right],\\
&p_{2k}^{+}:=\E \left[\psi_{2k}(M_k-\bar m-V(0) )\right],\quad p_{2k}^{-}:=\E \left[\psi_{2k+1}(M_k-\bar m-V(0) )\right].
\end{align*}
In light of \eqref{explicit-small}, for $1\leq k \leq L$, $p_{2k}< p_{2k-1}^{-}<p_{2k-1}<p_{2k-1}^{+}<p_{2k-2}$, and $p_{2k+1}< p_{2k}^{-}<p_{2k}<p_{2k}^{+}< p_{2k-1}$.

The following lemma says that  $(H, V)$ is regularly homogenizable under assumption \eqref{explicit-small}, that is, when the oscillation of $V$ is smaller than the depth of any well in the graph of $H$.

\begin{lemma}\label{localization}
We have $(H,V)$ is regularly homogenizable and the formula of $\overline H$ is given as follows. 

\noindent {\rm (1)} For $p\in [p_{2k}^{+}, p_{2k-2}^{-}]$ where $1\leq k \leq L$, $\overline H(p)$ is given by
\begin{equation}\label{local1}
\begin{cases}
p=\E \left[\psi_{2k-1}(\overline H(p)-V(0) )\right]  \quad \text{for $p\in [p_{2k-1}^{+},  p_{2k-2}^{-}]$}\\
\overline H(p)\equiv m_k    \quad \text{for $p\in [p_{2k-1}^{-}, p_{2k-1}^{+}]$}\\
p=\E \left[\psi_{2k}(\overline H(p)-V(0)) \right]  \quad \text{for $p\in [p_{2k}^{+},  p_{2k-1}^{-}]$}.
\end{cases}
\end{equation}
If $k=1$, the first equality becomes $p=\E \left[\psi_{1}(\overline H(p)-V(0) \right]$ for  $p\in [p_{1}^{+}, \infty)$.

\noindent{\rm(2)} For $p\in [p_{2k+1}^{+}, p_{2k-1}^{-}]$ where $1\leq k \leq L$, $\overline H(p)$ is given by 
\begin{equation}\label{local2}
\begin{cases}
p=\E \left[\psi_{2k}(\overline  H(p)-V(0) )\right]  \quad \text{for $p\in [ p_{2k}^{+},p_{2k-1}^{-} ]$}\\
\overline H(p)\equiv M_k-\bar m    \quad \text{for $p\in [p_{2k}^{-}, p_{2k}^{+}$}]\\
p=\E \left[\psi_{2k+1}(\overline H(p)-V(0) )\right]  \quad \text{for $p\in [ p_{2k+1}^{+},p_{2k}^{-} ]$}.
\end{cases}
\end{equation}

\noindent{\rm (3)}  For $p\leq \E \left[\psi_{2L+1}(-V(0, ) \right]$,  $\overline H(p)$ is given by 
$$
\begin{cases}
\overline H(p)\equiv  0  \quad \text{ for $p\in  [\E \left[\Psi(-V(0)) \right], \   \E \left[\psi_{2L+1}(-V(0)) \right]]$}\\
p=\E \left[\Psi(\overline H(p)-V(0)) \right] \quad \text{for $p \leq \E \left[\Psi(-V(0) \right]$}.
\end{cases}.
$$

\end{lemma}

\begin{proof}  We only prove (1) as the proofs of (2) and (3)  are similar.   It suffices to verify  the middle equality in  (\ref{local1}).  Other two equalities are obvious due to the existence of  sublinear solutions to the cell problem.   Our goal is to show that  for $p\in  [p_{2k-1}^{-}, p_{2k-1}^{+}]$
\begin{equation}\label{e11}
\P \left[\lim_{\lambda\to 0}|\lambda v_{\lambda}(0, p)+m_k|=0\right]=1,
\end{equation}
where  $v_{\lambda}(\cdot, p)\in C^{0,1}(\R)$ is the solution of \eqref{Cell-p}.

Let $\tilde H\in C^{0,1}(\R)$ be a  function satisfying that  $\tilde H\geq H$ (see the figure below), and

\[
\begin{cases}
\tilde H=H \quad \text{on}\ [p_{2k}, p_{2k-2}],\\
\tilde H \quad \text{is strictly increasing on} \ [p_{2k-1}, \infty),\\
\tilde H \quad  \text{is strictly decreasing on}\ (-\infty,  p_{2k-1}].
\end{cases}
\]
\begin{center}
\includegraphics[scale=0.7]{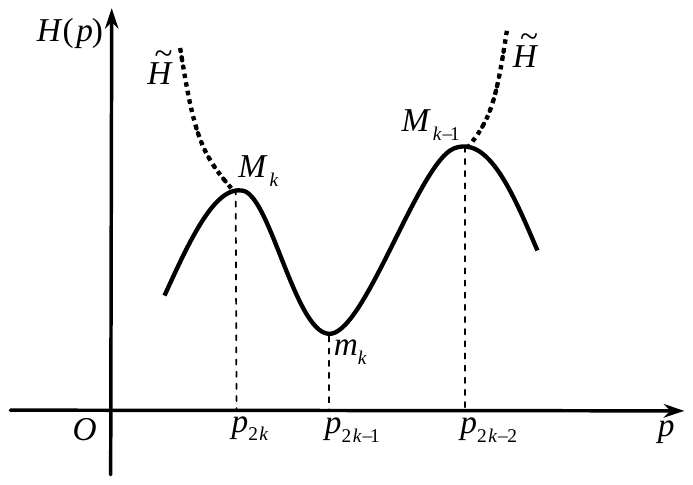}
\captionof{figure}{Construction of $\tilde H$}.
\end{center}

Owing to the previous lemma,   $(\tilde H, V)$ is regularly homogenizable and
$$
\overline {\tilde H}(p)=m_k  \quad \text{for $p\in   [p_{2k-1}^{-}, p_{2k-1}^{+}]$}.
$$
Now fix $p\in [p_{2k-1}^{-}, p_{2k-1}^{+}]$. Thus, for any $R>0$,
\begin{equation}\label{erg1}
\P \left[ \limsup_{\lambda\to 0}\max_{|y|\leq {R/\lambda}}|\lambda \tilde v_{\lambda}(y, p)+m_k|=0 \right] = 1,
\end{equation}
where $\tilde {v}_{\lambda}(\cdot, p)\in  C^{0,1}(\R)$ is the unique solution to
$$
\lambda {\tilde v}_{\lambda}+\tilde H(p+{\tilde v}_{\lambda}')+V(y)=0  \quad \text{in  $\R$}.
$$
It is a routine fact that
$$
\sup_{\R}|\lambda v_{\lambda}(\cdot, p)|,\quad  \sup_{\R}|\lambda \tilde v_{\lambda}(\cdot, p)|\leq H(p)+\bar m.
$$
Due to  (\ref{explicit-small}) and (\ref{erg1}), it is clear that, for fixed $V$ and $R>0$, there exists $\lambda(R, V)>0$ such that when $\lambda\leq \lambda(R, V)$,
$$
p+{\tilde v}_{\lambda}'(y, p)\in  (p_{2k},  p_{2k-2})   \quad \text{for $y\in   B_{R/\lambda}$}.
$$
So ${\tilde v}_{\lambda}(\cdot,  p)$ is also a viscosity solution to \eqref{Cell-p} in $B_{R/\lambda}$.
Hence according to Lemma~\ref{lem:comparison-R},  
$$
|\lambda {\tilde v}_{\lambda}(0,p)-\lambda v_{\lambda}(0,p)|\leq   {C\over R},
$$
where $C\geq 1$ depends only on $H$ and $\bar m$.    
This completes the proof of~(\ref{e11}).  
\end{proof}

%Since $p_{k}^{+}\leq p_{k-1}^{-}$,  as an immediate corollary,  we obtain that $(H,V)$ is regularly homogenizable for all $p\in  \R$  if the oscillation of $V$ is smaller than the depth of any well.  Denote
%$$
%D_{R}^{+}=\min_{1\leq i\leq L}\{M_i-m_i\},   \qquad  D_{L}^{+}=\min_{1\leq i\leq L-1}\{M_L,  M_i-m_{i+1}\}
%$$
%and
%$$
%D_{R}^{-}=\min_{1\leq i\leq {\tilde L}}\{{\tilde M}_i-{\tilde m}_i\},   \qquad  D_{L}^{-}=\min_{1\leq i\leq {\tilde L}-1}\{{\tilde M}_{\tilde L}, {\tilde  M}_i-\tilde m_{i+1}\}
%$$
%
%\begin{corollary}\label{smallosc}
%
% Suppose
%$$
%\bar m\leq\min \{D_{R}^{+}, \  D_{L}^{+}, \  D_{R}^{-}, \   D_{L}^{-}\}.
%$$  Then $(H,V)$ is regularly homogenizable for all $p\in  \R$.  When $p\geq 0$ and $1\leq k\leq L$,   $\overline H(p)$ is given by the following formula (the formula for $p\leq 0$ is similar):
%$$
%\begin{cases}
%p=\E \left[\psi_1(\overline H (p)-V(0) \right]   \quad \text{when $p\geq p_{1}^{+}$}\\
%\overline H(p)\equiv m_k   \quad \text{when $p\in  [p_{2k-1}^{-}, p_{2k-1}^{+}]$}\\
%p=\E \left[\psi_{2k}(\overline H (p)-V(0) \right]   \quad \text{when $p\in [p_{2k}^{+}, p_{2k-1}^{-}]$}\\
%\overline H (p)\equiv M_k-\bar m    \quad \text{when $p\in  [p_{2k}^{-}, p_{2k}^{+}]$}\\
%p=\E \left[\psi_{2k+1}(\overline H (p)-V(0) \right]   \quad \text{when $p\in [p_{2k+1}^{+}, p_{2k}^{-}]$}\\
%\overline H(p)\equiv 0  \quad \text{when $0\leq p\leq \E \left[\psi_{2L+1}(-V(0) \right]$}.
%\end{cases}
%$$
%
%\end{corollary}

\appendix
\section{Auxiliary lemmas}
\subsection{Some general results for viscosity solutions in $1$-dimensional space}
\begin{lemma}\label{A1} 
Assume that $H\in C(\R)$ is coercive and $\min_{\R}H=H(0)=0$.     
For any $\mu\geq 0$, there exists a Lipschitz continuous viscosity solution $u$ to
\[
H(u')+V(y)=\mu  \qquad \text{in} \ \R
\]
such that
\[
u'\geq  0  \qquad \text{ for a.e.} \ y \in  \R.
\]
\end{lemma}

\begin{proof}   We present the proof in two steps.

{\it  Step 1.}   We first assume that $V$ is periodic with period $1$ and $\mu>0$.  
Let $\overline H$ be the corresponding effective Hamiltonian.   
It is easy to see that $\overline H(0)=0$.   Choose $p_{\mu}>0$ such that $\overline H(p_{\mu})=\mu>0$.
Let $v \in C^{0,1}(\R)$ be a periodic viscosity solution to the cell problem
\[
H(p_{\mu}+v')+V(y)=\mu  \quad \text{in} \ \R.
\]
We claim that $u=p_{\mu}y+v$ satisfies that
\begin{equation}\label{A1-1}
u'>0 \qquad \text{for a.e.} \ y \in \R.
\end{equation}
Assume not,  then there exists $x_1\in \R$ such that $u'(x_1)\leq 0$.   
Since 
\[
p_{\mu}=\int_{x_1}^{x_1+1}u'(x)\,dx>0,
\]
there exists $x_2>x_1$ such that   $u'(x_2)>0$.  
Due to Lemma \ref{lem:meanvalue},    we may find $x_3 \in [x_1, x_2)$ such that $0 \in D^- u(x_3)$. By definition of viscosity solutions,
$$
H(0)+V(x_3)\geq \mu>0,
$$
which is absurd. Thus \eqref{A1-1} holds.

{\it Step 2:}  Now for $n\in \N$,  let $V_n\in C(\R)$ satisfy that
\begin{itemize}
\item $V_n(y)=V(y)$ for $|y|\leq n$.
\item  $V_n(y+2n)=V_n(y)$ for all $y\in  \R$,  $\max_{\R}V_n=0$ and $\max_{\R}|V_n|\leq \sup_{\R}|V|$;
\end{itemize}
Then owing to Step 1,  for $\mu\geq  0$ and $n\in \N$,  there exists   $u_n\in  C^{0,1}(\R)$  such that
\[
H(u_n')+V_n(y)=\mu+\frac{1}{n}  \quad \text{in} \ \R,
\]
and $u_n' >0$ a.e. in $\R$.

Due to the  coercivity of $H$ and the uniform boundedness of $\{V_n\}$,   $u_{n}$ is equi-Lipschitz continuous in $\R$.   Without loss of generality,  we may assume that
$$
u_n\to u   \quad \text{locally uniformly in} \  \R.
$$
By usual stability results of viscosity solutions, $u$ satisfies all the requirements of the lemma. 
\end{proof}

\begin{lemma}\label{A2}  
Assume that $H$ satisfies {\rm (H1)-(H2)} and levels set of $V$ have no cluster points.   Let $u\in C^{0,1}(\R)$ be a  viscosity solution of
\[
H(u')+V(y)=\mu\geq 0  \quad \text{in} \ \R
\]
and $u'\geq 0$ a.e. in $\R$. 
Then there  exists a strictly increasing sequence $\{b_i\}_{i\in  \Z}$ such that  $\lim_{i\to \pm \infty}b_i=\pm \infty$ and for $I_i:=(b_i, b_{i+1})$,  $u\in C^1(I_i)$ and
$$
u'|_{I_i}=\psi_{k_i}(\mu-V)  \quad \text{for some  $k_i\in \{1,2,\ldots, 2L+1\}$}.
$$
\end{lemma}

\begin{proof}
We claim that for any $y \in  \R$, there exists $\delta_y>0$ and $l_y,r_y\in  \{1,2,\ldots, 2L+1\}$ such that
$$
u'=
\begin{cases}
\psi_{r_y}(\mu-V)   \quad \text{in   $(y, y+\delta_y)$}\\
\psi_{l_y}(\mu-V(y)) \quad \text{in   $(y-\delta_y,y)$}
\end{cases}
$$
Let us prove the first equality.  Assume by contradiction that  there exist a decreasing sequence $\{y_n\}$ converging to $y$ and two numbers $k, k' \in \{1,2,\ldots,2L+1\}$ such that  $k>k'$, and for all $n\in \N$,
\begin{equation}\label{A2-1}
\begin{cases}
u'(y_{2n-1})=\psi_{k}(\mu-V(y_{2n-1}))\in [p_{k}, p_{k-1}], \\
 u'(y_{2n})=\psi_{k'}(\mu-V(y_{2n}))\in [p_{k'}, p_{k'-1}].
  \end{cases}
\end{equation}
This together with Lemma \ref{lem:meanvalue}  yield the existence of a sequence $\{z_n\}$ such that $z_n \in [y_{n+1},y_n]$ with $p_{k-1} \in D^+ u(z_{2n-1})$, and $p_{k-1} \in D^- u(z_{2n})$ for all $n\in \N$. Hence
\begin{equation}\label{A2-2}
H(p_{k-1})+V(z_{2n-1}) \leq \mu \leq H(p_{k-1})+V(z_{2n}).
\end{equation}
By the usual mean value theorem, there exists a further sequence $\{\bar z_n\}$ with $\bar z_n \in [z_{n+1},z_n]$ for all $n\in\N$ and
\begin{equation}\label{A2-3}
H(p_{k-1})+V(\bar z_n)=\mu,
\end{equation}
which implies that $y$ is a cluster point of $V$, and hence, contradiction. Therefore, \eqref{A2-1} holds, and furthermore $l_y,r_y$ are unique. Set
\[
A=\{y\in \R\,:\,l_y \neq r_y\}.
\]
By the same reason like the above step, $A$ has no cluster points and we can find a strictly increasing sequence $\{b_i\}_{i\in \Z}$ such that $\lim_{i\to \pm \infty}b_i=\pm \infty$ and $A \subseteq \{b_i\}_{i\in \Z}$.
\end{proof}

\subsection{Homotopy between solutions}
Take $f\in \A(H,V,\mu)$ and $b_1<b_2<b_3$ such that for $i=1,2$
\begin{equation}\label{A3-1}
f|_{(b_i,b_{i+1})}=\psi_{k_i}(\mu-V) \quad \text{for some} \ k_i \in \{1,2,\ldots,2L+1\}.
\end{equation}
Denote $k:=\min\{k_1,k_2\}$ and
\[
\tilde f:=\begin{cases}
f \qquad &\text{in} \ \R \setminus (b_1,b_3),\\
\psi_k(\mu-V)\qquad &\text{in} \  (b_1,b_3).
\end{cases}
\]
\begin{lemma}\label{A3}
If
\begin{equation}\label{A3-2}
\{\mu-V(y)\,:\,b_1 <y<b_3\} \cap \{M_i,m_j\,:\, 1\leq i,j\leq L\}=\emptyset,
\end{equation}
then $\tilde f \in \A(H,V,\mu)$.
\end{lemma}

\begin{proof}
Assume $k_1<k_2$.  Due to \eqref{A3-2},   both $\psi_{k_1}(\mu-V)$ and $\psi_{k_2}(\mu -V)$ are well defined in  $(b_1, b_3)$.   Let
$$
q_1:=\psi_{k_1}(\mu-V(b_2))    \quad \mathrm{and}  \quad  q_2:=\psi_{k_2}(\mu-V(b_2)).
$$
Clearly   $q_1>q_2$, $D^+u(b_2)=[q_2,q_1]$, and in light of \eqref{A3-2} for any $p\in (q_1,  q_2)$,
$$
H(p)<\mu-V(b_2).
$$
We actually can infer furthermore that $\psi_{k_1}$ is strictly increasing, and $\psi_{k_2}$ is strictly decreasing. For any $y\in (b_2,b_3)$, set
$$
q_{1,y}:=\psi_{k_1}(\mu-V(y))    \quad \mathrm{and}  \quad  q_{2,y}:=\psi_{k_2}(\mu-V(y)).
$$
Duet to continuity and 2d topology (see the figure below),  one still has
\begin{equation}\label{A3-3}
\max_{p\in [q_{2,y},q_{1,y}]} H(p) = \mu- V(y),
\end{equation}
which yields that $\tilde f_y \in \A(H,V,\mu)$, where
\[
\tilde f_y = \begin{cases}
f \qquad &\text{in} \ \R\setminus (b_1,b_3),\\
\psi_{k_1}(\mu-V) \qquad &\text{in} \ (b_1,y),\\
\psi_{k_2}(\mu-V) \qquad &\text{in} \ (y,b_3).
\end{cases}
\]
Letting $y \to b_3$ yields the desired result.
The proof for the case $k_1>k_2$ is similar hence omitted.
\end{proof}

\begin{center}
\includegraphics[scale=0.7]{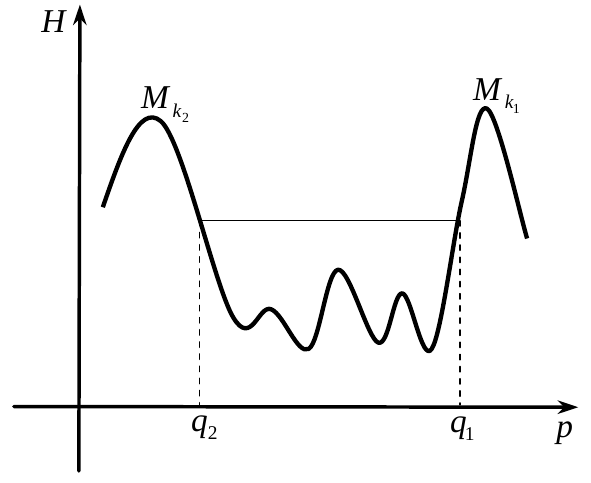}
\captionof{figure}{Position of $q_1$ and $q_2$}.
\end{center}

Take $f_1,f_2 \in \A(H,V,\mu)$. Assume there exist $a,b\in \R$ with $a<b$ such that 
\[
f_1 \geq f_2 \quad \text{in} \ I:=(a,b),\quad f_1=f_2 \quad \text{on} \ \R \setminus I.
\]
Pick $u_i \in C^{0,1}(\R)$ such that $u_i'=f_i$ and $u_i(a)=0$ for $i=1,2$.
It is straightforward that
\[
u_2(y) \leq u_1(y) \leq u_2(y)-u_2(b)+u_1(b) \quad \text{for} \ y \in [a,b].
\]
For any $c\in  [u_2(b),  u_1(b)]$ and $y\in I$, denote
$$
u_{c,*}(y):=\max\{u_2(y), \  u_1(y)-u_1(b)+c\},  \quad u_{c}^{*}(y):=\min\{u_1(y), \  u_2(y)-u_2(b)+c\}.
$$
Then $u_{c}^{*}\geq u_{c,*}$ in $I$, and  $u_{c}^{*}$ ($u_{c,*}$) are viscosity supersolution (subsolution) to \eqref{metric-mu} subject to
$$
u_{c}^{*}(a)=u_{c,*}(a)=0   \quad \mathrm{and} \quad u_{c}^{*}(b)=u_{c,*}(b)=c.
$$
For $y\in I$, define
\begin{equation*}
u_c(y):=\sup\left\{w(y)\,:\,w \ \text{is a subsolution of \eqref{metric-mu} and} \ u_{c,*}\leq w \leq u_c^* \ \text{in} \ I\right\}.
\end{equation*}
Also set $f_c=f_c(f_1,f_2,I)$ such that $f_c:=u_c'$ in $I$.

By abuse of notation, we extend $f_c$ to the whole $\R$ as
\begin{equation*}
f_c=\begin{cases}
u_c' \qquad &\text{in} \ I,\\
f_1 \qquad &\text{on} \ \R \setminus I.
\end{cases}
\end{equation*}
\begin{lemma}\label{A4}
For any $c\in  [u_2(b),  u_1(b)]$, $f_c \in \A(H,V,\mu)$.
\end{lemma}

\begin{proof}
Let $\tilde u_c$ be the extension of $u_c$ to $\R$ as
\begin{equation}\label{A4-1}
\tilde u_c:=\begin{cases}
u_1 \qquad &\text{on} \ (-\infty,a],\\
u_c \qquad &\text{in} \ I,\\
u_1 - u_1(b)+c \qquad &\text{on} \ [b,\infty).
\end{cases}
\end{equation}
We now show that $\tilde u_c$ is a viscosity solution of \eqref{metric-mu}.
It is enough to check the definition of viscosity solutions at $y=a$ and $y=b$.
At $y=a$, we have that
\begin{equation}\label{A4-2}
\begin{cases}
u_1 \geq \tilde u_c \geq u_2 \qquad &\text{in} \ I,\\
u_1=u_2=\tilde u_c \qquad &\text{on} \ (-\infty,a],
\end{cases}
\end{equation}
and hence $D^- \tilde u_c(a) \subset D^- u_1(a)$, $D^+ \tilde u_c(a) \subset D^+ u_2(a)$. These of course imply that $\tilde u_c$ is a viscosity solution of \eqref{metric-mu} at $y=a$.

At $y=b$, it is also clear that
\begin{equation}\label{A4-3}
\begin{cases}
u_1 - u_1(b) \leq \tilde u_c  - \tilde u_c(b) \leq u_2 - u_2(b) \qquad &\text{in} \ I,\\
u_1-u_1(b)=u_2-u_2(b)=\tilde u_c - \tilde u_c(b) \qquad &\text{on} \ [b,\infty),
\end{cases}
\end{equation}
which gives that $\tilde u_c$ is a viscosity solution of \eqref{metric-mu} at $y=b$ by a similar argument like the above.
\end{proof}

\subsection{Approximation of potential $V$}
For $\ep>0$,  consider the approximation of $V$ by analytic functions:
$$
V_\ep(y)=  \frac{1}{\sqrt {2\pi \ep}}\int_{\R}e^{\frac{-(z-y)^2}{\ep}}V(z)\,dz.
$$
It is easy to check that $V_{\ep}:\R \to \R$ is also stationary.  

\begin{lemma}\label{A5}
The followings hold
\begin{itemize}
\item[(i)] $\lim_{\ep \to 0} \|V_\ep-V\|_{L^\infty(\R \times \Omega)} =0$.
\item[(ii)] The level sets of $V_\ep$ have no cluster points.
\end{itemize}
\end{lemma}

\begin{proof}
The first assertion is obvious.
As for (ii), if it were wrong, there would exist $y_0 \in \R$ such that $V_\ep^{(k)}(y_0)=0$ for all $k \in \N$.
Assume without loss of generality that $y_0=0$ and $V_\ep(0)=0$. Then
$$
\int_{\R}y^k e^{-\frac{y^2}{\ep}}V(y)\,dy=0   \quad \text{for all} \ k\geq 0.
$$
Using Fourier transform,  we obtain that $V\equiv 0$, which is absurd.
\end{proof}

\begin{lemma}\label{A6}  

\begin{multline*}
\P \big[ \mbox{for every unbounded interval $I\subset \R$}, \\ \  (\inf V,  \sup V)\subseteq   V(I):=\{V(y)\,:\,  y\in I\} \big]  = 1.
\end{multline*}
\end{lemma}
\begin{proof}
Using rational numbers,  it suffices to show that for any $c\in  (\inf V,  \sup V) \cap \Q$,  $a\in  \Q$, $I_{a}^{+}:=(a,\infty)$ and $I_{a}^{-}:=(-\infty,  a)$
$$
\P \left[  [c,  \sup V)\cap V(I_{a}^{+})=\emptyset \right]=\P \left[ (\inf V,  c]\cap V(I_{a}^{+})=\emptyset \right]=0
$$
and
$$
\P \left[  [c,  \sup V)\cap V(I_{a}^{-})=\emptyset \right]=\P \left[ (\inf V,  c]\cap V(I_{a}^{-})=\emptyset \right]=0.
$$
Let $g:=\mathbf{1}_{[c,\sup V)}$ and observe, by the ergodic theorem, that
\begin{equation*}
\P \left[ \lim_{L\to +\infty}  \frac{1}{L-a} \int_{a}^{L}  g(V(y))\,dy=\mathbb {E}( g(V(0)))>0 \right] =1.
\end{equation*}
This shows that $\P \left[  [c,  \sup V)\cap V(I_{a}^{+})=\emptyset \right]=0$. The proofs for the other equalities are similar. 
\end{proof}

\smallskip

\noindent{\bf Acknowledgements.}
We thank Yu-yu Liu for his tremendous help in drawing the pictures in this paper.
The second author is supported in part by NSF grant DMS-1361236.
The third author is supported in part by NSF CAREER award \#1151919.

\bibliographystyle{plain}
\bibliography{nonconvex}

\begin{thebibliography}{10}

\bibitem{AC}
S.~N. Armstrong and P.~Cardaliaguet.
\newblock {Quantitative stochastic homogenization of viscous Hamilton-Jacobi
  equations}.
\newblock {\em Comm. Partial Differential Equations}, to appear.

\bibitem{ACS}
S.~N. Armstrong, P.~Cardaliaguet, and P.~E. Souganidis.
\newblock {Error estimates and convergence rates for the stochastic
  homogenization of Hamilton-Jacobi equations}.
\newblock {\em J. Amer. Math. Soc.}, 27:479--540, 2014.

\bibitem{AS2}
S.~N. Armstrong and P.~E. Souganidis.
\newblock Stochastic homogenization of level-set convex {H}amilton-{J}acobi
  equations.
\newblock {\em Int. Math. Res. Not.}, 2013(15):3420--3449, 2013.

\bibitem{AT1}
S.~N. Armstrong and H.~V. Tran.
\newblock {Stochastic homogenization of viscous Hamilton-Jacobi equations and
  applications}.
\newblock {\em Analysis \& PDE}, to appear.

\bibitem{ATY2013}
S.~N. Armstrong, H.~V. Tran, and Y.~Yu.
\newblock {Stochastic homogenization of a nonconvex Hamilton-Jacobi equation},
  preprint, arXiv:1311.2029 [math.AP].

\bibitem{DS1}
A.~Davini and A.~Siconolfi.
\newblock Exact and approximate correctors for stochastic {H}amiltonians: the
  1-dimensional case.
\newblock {\em Math. Ann.}, 345(4):749--782, 2009.

\bibitem{EBook}
L.~C. Evans.
\newblock {\em Partial differential equations}, volume~19 of {\em Graduate
  Studies in Mathematics}.
\newblock American Mathematical Society, Providence, RI, 1998.

\bibitem{KRV}
E.~Kosygina, F.~Rezakhanlou, and S.~R.~S. Varadhan.
\newblock Stochastic homogenization of {H}amilton-{J}acobi-{B}ellman equations.
\newblock {\em Comm. Pure Appl. Math.}, 59(10):1489--1521, 2006.

\bibitem{KV}
E.~Kosygina and S.~R.~S. Varadhan.
\newblock Homogenization of {H}amilton-{J}acobi-{B}ellman equations with
  respect to time-space shifts in a stationary ergodic medium.
\newblock {\em Comm. Pure Appl. Math.}, 61(6):816--847, 2008.

\bibitem{LS1}
P.-L. Lions and P.~E. Souganidis.
\newblock Correctors for the homogenization of {H}amilton-{J}acobi equations in
  the stationary ergodic setting.
\newblock {\em Comm. Pure Appl. Math.}, 56(10):1501--1524, 2003.

\bibitem{LS2}
P.-L. Lions and P.~E. Souganidis.
\newblock Homogenization of ``viscous'' {H}amilton-{J}acobi equations in
  stationary ergodic media.
\newblock {\em Comm. Partial Differential Equations}, 30(1-3):335--375, 2005.

\bibitem{MN}
I.~Matic and J.~Nolen.
\newblock {A sublinear variance bound for solutions of a random
  Hamilton--Jacobi equation}.
\newblock {\em Journal of Statistical Physics}, 149(2):342--361, 2012.

\bibitem{Q}
J.~L. Qian.
\newblock {Two approximations for effective Hamiltonians arising from
  homogenization of Hamilton-Jacobi equations}, UCLA CAM reports (2003),
  03--39.

\bibitem{RT}
F.~Rezakhanlou and J.~E. Tarver.
\newblock Homogenization for stochastic {H}amilton-{J}acobi equations.
\newblock {\em Arch. Ration. Mech. Anal.}, 151(4):277--309, 2000.

\bibitem{Sch}
R.~W. Schwab.
\newblock Stochastic homogenization of {H}amilton-{J}acobi equations in
  stationary ergodic spatio-temporal media.
\newblock {\em Indiana Univ. Math. J.}, 58(2):537--581, 2009.

\bibitem{S}
P.~E. Souganidis.
\newblock Stochastic homogenization of {H}amilton-{J}acobi equations and some
  applications.
\newblock {\em Asymptot. Anal.}, 20(1):1--11, 1999.

\end{thebibliography}

\end{document}